\documentclass[11pt]{article}
\usepackage[american]{babel}
\usepackage{microtype}
\usepackage{graphicx}
\usepackage{pdfsync}
\usepackage{amsmath,amssymb,amsthm,amsfonts,mathrsfs}
\usepackage{amssymb}
\usepackage{verbatim}
\usepackage{setspace}
\usepackage{hyperref}
\usepackage[a4paper]{geometry}

\newtheorem{thm}{Theorem}[section]

\newtheorem{lem}[thm]{Lemma}

\theoremstyle{definition}

\theoremstyle{remark}

\numberwithin{equation}{section}

\DeclareMathSymbol{\C}{\mathalpha}{AMSb}{"43}

\newcommand{\eps}{\epsilon}

\newcommand{\R}{{\mathbb{R}}}
\newcommand{\h}{{\mathcal{H}}}

\newcommand{\dt}{\,\mathrm{d}t}
\newcommand{\dx}{\,\mathrm{d}x}

\newcommand{\ds}{\,\mathrm{d}S}

\allowdisplaybreaks

\newcommand{\bsub}{\begin{subequations}}
\newcommand{\esub}{\end{subequations}$\!$}

\begin{document}
\title{Mass Concentration and Local Uniqueness of Ground States for $L^2$-subcritical Nonlinear Schr\"{o}dinger Equations\thanks{This work is partially supported by NSFC under Grant No. 11671394. }}
\author{
Shuai Li$^{a\ b}$
\ and\   Xincai Zhu$^{a\ b}$
\thanks{Corresponding author.
E-mail: lishuai\_wipm@outlook.com (S. Li); zhuxc68@163.com (X. Zhu).}
\\
\small \it $^{a}$ Wuhan Institute of Physics and Mathematics, Chinese Academy of Sciences,\\
\small \it Wuhan 430071, P. R. China\\
\small \it$^{b}$ University of Chinese Academy of Sciences, Beijing 100049, P. R. China\\
}

\date{\today}
\smallbreak \maketitle

\begin{abstract}
We consider ground states of  $L^2$-subcritical  nonlinear Schr\"{o}dinger equation \eqref{equ:u},
which can be described equivalently by  minimizers of the following  constraint minimization problem
\begin{equation*}\label{def:e}
  e(\rho):=\inf\big\{E_{\rho}(u):\,u\in \h(\R^d),\|u\|_2^2=1\big\}.
\end{equation*}
The  energy functional $E_{\rho}(u)$ is defined by
\begin{equation*}\label{def:E}
  E_{\rho}(u):=\frac{1}{2}\int_{\mathbb{R}^d}|\nabla u|^2\dx
  +\frac{1}{2}\int_{\mathbb{R}^d}V(x)|u|^2\dx
  -\frac{\rho^{p-1}}{p+1}\int_{\mathbb{R}^d}|u|^{p+1}\dx,
\end{equation*}
where  $d\geq1$, $\rho>0$, $p\in\big(1, 1+\frac{4}{d}\big)$ and $0\leq V(x)\to\infty$ as $|x| \to\infty$.
We present a detailed analysis on the concentration behavior of  ground states as $\rho\to\infty$, which extends the concentration results shown in \cite{M}. %
Moreover, the uniqueness of nonnegative ground states is also proved when $\rho$ is large enough.
\end{abstract}



\section{Introduction}
In this paper, we study the following time-independent  nonlinear Schr\"{o}dinger equation
\begin{equation}\label{equ:u}
  -\Delta u+V(x)u=\mu u+\rho^{p-1}u^p,\,\,\,u\in \h(\R^d),
\end{equation}
where $d\geq1$, $\mu\in\R$, $p\in\big(1, 1+\frac{4}{d}\big)$,
and $\rho>0$ describes  the strength of the attractive interactions.

The space $\h(\R^d)$ is defined as
\begin{equation}\label{def:H}
   \h(\R^d) := \Big \{u(x)\in H^1(\R^d)\,\  \Big|\,\  \int _{\R^d}  V(x)|u(x)|^2 \dx<\infty \Big\},
\end{equation}
and the associated norm is given by
$\|u\|_{\h}=\big\{\int_{\R^d}\big(|\nabla u(x)|^2+[1+V(x)]|u(x)|^2\big)\dx\big\}^\frac{1}{2}$.
Equation \eqref{equ:u}  arises in Bose-Einstein condensates (BEC) and nonlinear optics.
Especially, when $p=3$ and $d=1$, it is the well known time-independent Gross-Pitaevskii (GP) equation which describes the one-dimensional BEC problem, see, e.g., \cite{Gross1,Gross2, Pitaevskii} and the references therein.
From the physical point of view, we assume that the trapping potential $V(x)\geq0$ satisfies
\begin{equation}\label{con:V1}
\text{$V(x)\in L_{loc}^\infty(\R^d)\cap C^\alpha(\R^d)$ with $\alpha\in(0,1)$, $\inf\limits_{x\in\R^d}V(x)=0$ and $\lim\limits_{|x|\to\infty}V(x) = \infty$. }
\end{equation}

It is well known that,
 a minimizer of the following  minimization problem solves equation \eqref{equ:u} for some  suitable Lagrange multiplier $\mu$,
\begin{equation}\label{def:e}
  e(\rho):=\inf\big\{E_{\rho}(u):\,u\in \h(\R^d),\|u\|_2^2=1\big\},
\end{equation}
where $E_{\rho}(u)$ is the  Gross-Pitaevskii (GP) energy functional defined by
\begin{equation}\label{def:E}
  E_{\rho}(u):=\frac{1}{2}\int_{\mathbb{R}^d}|\nabla u|^2\dx
  +\frac{1}{2}\int_{\mathbb{R}^d}V(x)|u|^2\dx
  -\frac{\rho^{p-1}}{p+1}\int_{\mathbb{R}^d}|u|^{p+1}\dx.
\end{equation}
Equivalence between ground states of equation \eqref{equ:u} and constraint minimizers of \eqref{def:e} is proved in Theorem \ref{Thm:g.m}.
To discuss equivalently ground states of \eqref{equ:u}, in this paper,
we shall therefore focus on investigating  \eqref{def:e}, instead of  \eqref{equ:u}.
On the other hand, as for the general constraint $\|u\|_2=N\in(0,\infty)$, one can check that this latter case can be easily reduced to \eqref{def:e}, by
 minimizing \eqref{def:E} under the constraint $\|u\|_2=1$ but simply replacing $\rho$ by $\frac{\rho}{N}$.

When $p>1+\frac{4}{d}$, \eqref{def:e} is the so called $L^2$-supcritical problem (also known as mass supcritical problem).
Taking a suitable trial function and substituting it into \eqref{def:E}, one can check that problem \eqref{def:e} admits no minimizer for any $\rho\in(0,\infty)$  under this case.
For the case $p=1+\frac{4}{d}$, \eqref{def:e} is known as the $L^2$-critical problem (also called mass critical problem).
Recently, some interesting results  on this $L^2$-critical problem  were obtained by Guo and his co-authors (cf. \cite{GLW}-\cite{GZZ}).
Roughly speaking, the authors proved in \cite{GS} that there exists a finite value $\rho^*$ such that \eqref{def:e} admits minimizers if and only if $\rho<\rho^*$ (see also \cite{Bao2} for similar results).
The threshold value $\rho^*$ is determined by $\|w\|_2$, where $w$ is the unique (up to translations) positive radially symmetric solution of the following nonlinear scalar field equation (cf.  \cite{GNN,Kwong,LN,mcleod})
\begin{equation}\label{equ:w}
\Delta w-w+w^p=0,\,\,\  w \in H^1(\R^d).
\end{equation}
The concentration behavior of minimizers as $\rho \nearrow \rho^*$ was also analyzed in \cite{GS, GWZZ, GZZ} under different types of trapping potentials.
Furthermore, the local uniqueness of minimizers as $\rho\nearrow \rho^*$ was proved in \cite{GLW}, where the trapping potential is a class of homogeneous functions.

As for the $L^2$-subcritical case, $i.e.$, $1<p<1+\frac{4}{d}$, problem \eqref{def:e} admits minimizers for any $\rho\in(0,\infty)$, see, e.g., \cite{C,GZZJDE,M,Rabinowitz,Zeng2}.
Some qualitative properties of minimizers for \eqref{def:e}, such as uniqueness, concentration behavior and symmetry, were also studied in \cite{GZZJDE, M,Zeng2} and the references therein.
In detail, M. Maeda showed in \cite{M} that  minimizers of \eqref{def:e} are unique when $\rho$ is small enough
and the minimizers  must concentrate at a global minimum of $V(x)$ as $\rho\to\infty$.
Further,   Guo, Zeng and Zhou  presented in \cite{GZZJDE} a detailed analysis on the concentration behavior of minimizers for $e(\rho)$ with $d=2$ as $q\nearrow3$, and more recently, Zeng has generalized these results in \cite{Zeng2}.

Motivated by the works mentioned above,
in  this paper we focus on proving the local uniqueness of minimizers  for the $L^2$-subcritical problem \eqref{def:e} as $\rho\to\infty$.
Towards this purpose, it is necessary to analyze the concentration behavior of minimizers as $\rho\to\infty$.
However,  the concentration results shown in \cite{M} are not enough, and we therefore need to give a more detailed analysis on the limit behavior of minimizers for $e(\rho)$ as $\rho\to\infty$.
Besides, the equivalence between ground states of \eqref{equ:u} and constraint minimizers of \eqref{def:e} is also addressed.

Before stating our results, we need to introduce the following classical Gagliardo-Nirenberg type inequality  (cf. \cite{W})
\begin{equation}\label{ineq:GNw}
 C_{GN}
 \leq \frac{\|\nabla u\|_2^{\frac{d}{2}(p-1)}\|u\|_2^{p+1-\frac{d}{2}(p-1)}}{\|u\|_{p+1}^{p+1}} \,\,\,\text{for any}\,\,\,u\in H^1(\R^d)\setminus\{0\},
\end{equation}
where
\begin{equation}\label{ineq:GNw.C}
\begin{split}
C_{GN}
:
=&\|w\|_2^{p-1}\Big(1-\frac{p-1}{p+1}\frac{d}{2}\Big)\Big[\frac{d(p-1)}{2(p+1)-d(p-1)}\Big]^{\frac{d}{4}(p-1)},\\
\end{split}
\end{equation}
and $w$ is the unique positive solution of \eqref{equ:w}.
The equality in \eqref{ineq:GNw} is attained at $w$.
Applying the following Pohozaev identity of  \eqref{equ:w} (cf. \cite[Lemma 8.1.2]{C})
\begin{equation}\label{Pid:w}
(d-2)\int_{\R^d}|\nabla w|^2\dx+d\int_{\R^d}w^2\dx=\frac{2d}{p+1}\int_{\R^d}w^{p+1}\dx,
\end{equation}
one can deduce from \eqref{equ:w}  that $w$ satisfies
\begin{equation}\label{id:w}
\int_{\R^d}|\nabla w|^2\dx=\frac{d}{2}\frac{p-1}{p+1}\int_{\R^d}w^{p+1}\dx=\frac{d(p-1)}{2(p+1)-d(p-1)} \int_{\R^d}w^2\dx.
\end{equation}
Note also from \cite[Proposition 4.1]{GNN} that $w(x)$ decays exponentially in the sense that
 \begin{equation} \label{decay:w}
w(x) \, , \ |\nabla w(x)| = O(|x|^{-\frac{d-1}{2}}e^{-|x|}) \,\
\text{as} \,\ |x|\to \infty.
\end{equation}

Our first result is concerned with the equivalence between minimizers of \eqref{def:e} and ground states of \eqref{equ:u}.
For convenience, we introduce some notations in advance.
For any given  $\rho\in(0,\infty)$,
the set of nontrivial weak solutions for \eqref{equ:u} is defined by
\begin{equation*}
S_{\mu,\rho}:=\Big\{
  u\in \h\setminus \{(0)\}:
  \ \langle F_{\mu,\rho}'(u),\varphi\rangle=0,\, \forall \ \varphi\in \h
  \Big\},
\end{equation*}
where the energy  functional $F_{\mu,\rho}(u)$ is   defined  as
\begin{equation}\label{def:F}
\begin{split}
F_{\mu,\rho}(u)
:=&\frac{1}{2}\int_{\R^d}|\nabla u|^2\dx
  +\frac{1}{2}\int_{\R^d}\big(V(x)-\mu\big)|u|^2\dx
  -\frac{\rho^{p-1}}{p+1}\int_{\R^d}|u|^{p+1}\dx.\\
\end{split}
\end{equation}
Further, the set of ground states for \eqref{equ:u}  is given by
\begin{equation}
G_
{\mu,\rho}:=\Big\{
  u\in S_{\mu,\rho}:\,F_{\mu,\rho}(u)\leq F_{\mu,\rho}(v)\, \text{ for all }\, v \in S_{\mu,\rho}
  \Big\}.
\end{equation}
Moreover,  the set of minimizers for $e(\rho)$ is defined as
\begin{equation}\label{def:set.mini}
  M_\rho:=\Big\{ u_\rho\in\h(\R^d):\,\text{$u_\rho$ is a minimizer of $e(\rho)$}  \Big\}.
\end{equation}
Our first result is stated as the following theorem.

\begin{thm}\label{Thm:g.m}
Suppose $ V(x)$ satisfies \eqref{con:V1}.
Then we have follows.
\begin{itemize}
\item [(i).] For $a.e.$ $\rho\in(0,\infty)$, all minimizers of $e(\rho)$ satisfy equation \eqref{equ:u} with a fixed Lagrange multiplier $\mu=\mu_\rho$.
\item [(ii).] For $a.e.$ $\rho\in(0,\infty)$,  $G_{\mu_\rho,\rho}=M_\rho$.
\end{itemize}
\end{thm}

Theorem \ref{Thm:g.m} indicates that, for $a.e.$ $\rho\in(0,\infty)$, there exists a unique $\mu=\mu_\rho$ such that equation \eqref{equ:u} admits ground states, which are equivalent to minimizers of $e(\rho)$ in the sense that $G_{\mu_\rho,\rho}=M_\rho$.
Theorem \ref{Thm:g.m} is largely inspired by some similar conclusions on different types of problems, such as  \cite[Chapter 8]{C}, \cite[Theorem 1.1]{GWZZ} and \cite[Theorem 1.2]{LXZ}.
The proof of Theorem \ref{Thm:g.m} is left to Appendix \ref{sec:g.m}.

We next focus on analyzing the limit behavior of minimizers as $\rho\to\infty$.
Since $|\nabla |u||\leq |\nabla u|$ holds for $a.e.$ $x\in\R^d$,
without loss of generality, we always suppose minimizers of $e(\rho)$ are nonnegative.
Motivated by \cite{GS}-\cite{GZZ},
in order to analyze the blow-up behavior of minimizers as $\rho \to \infty$, some additional assumptions on  $V(x)$ are required.

\noindent {\bf Definition 1.1.} $h(x)$ in $\mathbb{R}^d$ is homogeneous of degree $q\in \mathbb{R}^{+}$ (about the origin), if there exists some $q>0$ such that
\[ h(tx)=t^{q}h(x)\,\,\,\hbox{in}\,\,\,\R^d~~\hbox{for any }t>0.\]
This definition indicates that if $h(x)\in C(\R^d)$ is homogeneous of degree $q>0$, then
\[\text{$0\le h(x)\le C|x|^{q}$ in $\R^d$, where $C:=\max\limits_{x\in\partial B_1(0)}h(x)$,} \]
because $h(\frac{x}{|x|})\le C$ for any $x\in\R^d\setminus{\{0\}}$.
Moreover, if $h(x) \to \infty$ as $|x|\to\infty $, then $0$ is the unique minimum point of $h(x)$.

Define the set of global minimum points of $V(x)$ by
\begin{equation}\label{def:beta.z}
  Z:=\big\{x\in\R^d:V(x)=0\big\}=\big\{x_1, x_2, \cdots , x_m\big\}, \,\ \text{where}\,\ m\geq1.
\end{equation}
We then assume that, $V(x)$ is almost  homogeneous of degree $r_i>0$ around each $x_i$.
Specifically, there exists some  $V_i(x)\in C^2_{\rm loc}(\R^d)$, which  is homogeneous of degree $r_i>0$ and satisfies $\lim\limits_{|x|\to\infty} V_i(x) = +\infty$, such that
\begin{equation}\label{con:V2}
\lim_{x\to0}\frac{V(x+x_i)}{ V_i(x)}=1,\,\,\,i=1,2,\cdots,m.
\end{equation}
Additionally, inspired by \cite{Grossi}, we define $Q_i(y) $ by
\begin{equation}\label{def:Q}
Q_i(y):=\int_{\R^d}V_i(x+y)w^2\dx,\,\,\,i=1,2,\cdots,m.
\end{equation}
Set
\begin{equation}\label{def:rho.p0}
  r:=\max_{1\leq i\leq m}r_i,\,\,\,\bar{Z}:=\big\{x_i\in Z: r_i= r\big\} \subset Z,
\end{equation}
and
\begin{equation}\label{def:rho.gamma}
\bar{\lambda}_0:=\min\limits_{ i\in \Gamma }\bar{\lambda}_i, \,\text{ where }\bar{\lambda}_i:=\min_{y\in \R^d}Q_i(y)\ \text{ and }\,\Gamma:=\big\{i: x_i\in \bar Z\big\}.
\end{equation}
Besides,  we also introduce some useful notations,
\begin{equation}\label{def:lam}
  \lambda:=\frac{1}{2}\frac{4-d(p-1)}{2(p+1)-d(p-1)},
\end{equation}
\begin{equation}\label{def:Q}
Q(y):=\int_{\R^d}V_0(x+y)w^2\dx,\,\,\,\text{where $V_0(x):=V_i(x)$ and $i$ satisfies $\bar{\lambda}_i=\bar{\lambda}_0$},
\end{equation}
and
\begin{equation}\label{def:rho.z0.y0}
Z_0:=\big\{x_i\in \bar{Z}: \bar{\lambda}_i= \bar{\lambda}_0\big\},\,\,\,K_0:=\big\{y: Q(y)=\bar{\lambda}_i=\bar{\lambda}_0\},
\end{equation}
where $Z_0$ denotes the set of the flattest global minimum points of $V(x)$.
Stimulated by \cite{Grossi2,GS,M,Wang}, we  now give the following theorem on the blow-up behavior of nonnegative minimizers as $\rho\to\infty$.

\begin{thm}\label{thm:uk}
Suppose $ V(x)\in C^2(\R^d)$ satisfies \eqref{con:V1} and \eqref{con:V2},
and there exists a constant $\kappa>0$ such that
\begin{equation}\label{con:V3}
 V(x)\leq C e^{\kappa |x|}\,  \text{  if }\,  |x| \text{ is large}.
\end{equation}
Set $a^*:=\|w\|_2^2$, where $w$ is the unique positive solution of \eqref{equ:w}.
Let $u_k$ be a nonnegative minimizer of $e(\rho_k)$, where $\rho_k\to \infty$ as $k\to\infty$.
Then there exists a subsequence, still denoted by $\{u_k\}$, such that $u_k$ satisfies
\begin{equation}\label{lim:u}
\bar{u}_k(x):=\sqrt{a^*}\varepsilon_k^\frac{d}{2}u_k(\varepsilon_kx+x_k) \to w(x)\,\,\,\text{uniformly in}\,\,\,\R^d\,\,\,\text{as}\,\,\,\,k\to\infty,
\end{equation}
where
\begin{equation}\label{def:eps}
\varepsilon_k:=\Big(\frac{\rho_k}{\sqrt{a^*}}\Big)^{-\frac{2(p-1)}{4-d(p-1)}},
\end{equation}
and $x_{k}$ is the unique local maximum point of $u_k$ satisfying
\begin{equation}\label{lim:x}
  \frac{x_k-x_0}{\varepsilon_k}\to y_0\,\,\,\text{ for some $x_0\in Z_0$ and $y_0\in K_0$}\,\,\,\,\,\text{as}\,\,\,\,k\to\infty.
\end{equation}
Furthermore, $u_k$ decays exponentially in the sense that
\begin{equation}\label{decay:u.exp}
\bar{u}_k(x)\le Ce^{-\frac{|x|}{2}} \,\ \text{and}\,\,\ |\nabla \bar{u}_k|\le Ce^{-\frac{|x|}{4}} \,\ \text{as}\,\,\ |x|\to\infty,
\end{equation}
where  $C>0$ is a constant independent of $k$.
\end{thm}

Theorem \ref{thm:uk} shows that minimizers of $e(\rho)$ must concentrate at one of the flattest global minimum points of $V(x)$ as $\rho\to\infty$.
Some results similar to \eqref{lim:u} were also obtained in \cite{M}. 
In Section \ref{sec:blowup}, we shall present a different proof for \eqref{lim:u} by employing the refined energy estimates in Lemma \ref{lem:est.e} and blow-up analysis in Lemma \ref{lem:wk}.
Here we point out that,  the $L^2$-subcritical nonlinearity term will lead to some difficulties on analyzing the limit behavior of  minimizers, due to that the Gagliardo-Nirenberg inequality cannot be used directly.
This is quite different from those obtained in \cite{GWZZ,GZZ}, where the $L^2$-critical problem is considered.
The key to solve these problems is to establish a refined energy estimate of $e(\rho)$ firstly, by which one can deduce the blow-up rates of $u_\rho$ as $\rho\to\infty$.
Towards this aim, we have to employ the fact that $e(\rho)\geq\tilde{e}(\rho)$, where $\tilde{e}(\rho)$ is a new minimization problem  defined in \eqref{def:tile}. 
Moreover, \eqref{lim:x} gives the convergence rate of the unique maximum point of each minimizer  as $\rho\to\infty $, which is based on a more precise energy estimate of $e(\rho)$.
In fact, we shall show that
\begin{equation}\label{est2:e}
  e(\rho)=-\lambda \Big(\frac{\rho}{\sqrt{a^*}}\Big)^\frac{4(p-1)}{4-d(p-1)}
  +\frac{\bar{\lambda}_0+o(1)}{a^*}\Big(\frac{\rho}{\sqrt{a^*}}\Big)^{-\frac{2r(p-1)}{4-d(p-1)}}\,\,\,\text{as $\rho\to\infty$,}
\end{equation}
where $a^*:=\|w\|_2^2$, $r$,  $\bar{\lambda}_0$ and $\lambda$ are respectively defined by \eqref{def:rho.p0}, \eqref{def:rho.gamma} and \eqref{def:lam}.

Motivated by the uniqueness results  addressed in \cite{Cao,Deng,Grossi,GLW}, we finally investigate the uniqueness of nonnegative minimizers for $e(\rho)$ as $\rho\to\infty$.
Towards this purpose, we require some additional conditions on $V(x)$.
Suppose  $V(x)$ admits a unique flattest minimum point $x_0$, $i.e.$,
\begin{equation}\label{con:V4}
\text{$Z_0$  contains only one element $x_0$, where $Z_0$ is defined in \eqref{def:rho.z0.y0}. }
\end{equation}
Further, we suppose that
\begin{equation}\label{asum:V.r2}
\text{$V(x)$ is homogeneous of degree $r\geq2$ near $x_0$. }
\end{equation}
Moreover, we also assume that there exists a constant $R_0$ small enough such that
\begin{equation}\label{con:V5}
\frac{\partial V(x+x_0)}{\partial x_i}=\frac{\partial  V_0(x)}{\partial x_i}+W_i(x)\,  \text{ and }\,\  |W_i(x)|\leq C|x|^{s_i}\,  \text{ in }\,    B_{R_0}(0),\,\
\end{equation}
where $x_0\in Z_0$, $V_0$ is given in \eqref{def:Q}, and $s_i>r-1$ for $i=1,2,\cdots,m$.
Under these assumptions,  our uniqueness results can be stated as the following theorem.

\begin{thm}\label{Thm:uniqueness} Suppose $V(x)\in C^2(\R^d)$ satisfies \eqref{con:V1}, \eqref{con:V2}, \eqref{con:V3} and \eqref{con:V4}-\eqref{con:V5}.  Moreover, we also assume that
\begin{equation}\label{con:V6}
\text{\text{$y_0$ is the unique and non-degenerate critical point of $Q(y)$},}
\end{equation}
where $Q(y)$ is defined by \eqref{def:Q}.
Then there exists a unique nonnegative minimizer for $e(\rho)$ as $\rho\to \infty$.
\end{thm}

Theorem \ref{Thm:uniqueness} indicates that problem \eqref{def:e} admits only one nonnegative minimizer when $\rho$ is large enough.
Together with Theorem \ref{Thm:g.m} and the uniqueness results given in \cite[Theorem 1.2]{M}, one can conclude that,
for any given $\rho$ where $\rho$ is small enough or large enough, there exists a unique $\mu=\mu_\rho$ such that
 \[\text{\eqref{equ:u} admits one and only one nonnegative ground state.}\]

Theorem \ref{Thm:uniqueness} is proved by establishing various types of local Pohozaev identities,
which is inspired by the \cite{Cao,Deng,Grossi,GLW}.
However,  the proof of Theorem \ref{Thm:uniqueness}  requires more involved and intricate calculations,  because of the general assumptions on dimension and trapping potentials.
Moreover, comparing with  discussing $L^2$-critical problem, the appearance of $L^2$-subcritical term also leads to some essential differences on deriving the second  Pohozaev identity.

This paper is organized as follows.
Section \ref{sec:blowup} is concerned with proving Theorem \ref{thm:uk}  on the limit behavior of minimizers for $e(\rho)$ as $\rho\to\infty$.
The main purpose of Section \ref{sec:unique} is to prove the local uniqueness of nonnegative minimizers by deriving local Pohozaev identities.
The proof of Theorem \ref{Thm:g.m} is left to Appendix \ref{sec:g.m}, and we also give some useful results on $\tilde{e}(\rho)$ in  Appendix \ref{sec:tile}.

\section{Mass concentration}\label{sec:blowup}
In this section, we shall prove Theorem \ref{thm:uk} on the limit behavior of minimizers for $e(\rho)$ as $\rho\to\infty$.
We  shall firstly establish the optimal energy estimates for $e(\rho)$,
and then present a detailed analysis on the limit behavior of minimizers as  $\rho\to\infty$.

\subsection{Refined energy estimates}\label{sec:enery.est}

The main purpose of this section is to establish the refined estimates of $e(\rho)$ by the following lemma.

\begin{lem}\label{lem:est.e}
Suppose $V(x)$ satisfies \eqref{con:V1}, and then we have
\begin{equation}\label{val:e}
 \lim_{\rho\to\infty} \frac{e(\rho)}{
 \Big(\frac{\rho}{\sqrt{a^*}}\Big)^\frac{4(p-1)}{4-d(p-1)}}
 =-\lambda,
\end{equation}
where $\lambda$ is given in \eqref{def:lam}, $a^*:=\|w\|_2^2$ and $w$ is the unique positive solution of \eqref{equ:w}.
\end{lem}
\noindent {\bf Proof.}
We start with the upper bound estimate  on the energy $e(\rho)$ as $\rho\to\infty$.
Suppose $\chi(x)\in C^\infty(\R^d)$ is a cut-off function satisfying $\chi(x)=1$ as $|x|\leq1$ and $\chi(x)=0$ as $|x|\geq2$.
Choose a trial function
\begin{equation}\label{def:utau}
u_{\tau}(x):=\frac{A_\tau}{\|w\|_{2}}\tau^\frac{d }{2}w(\tau x)\chi(x),
\end{equation}
where $\tau=\Big(\frac{\rho}{\sqrt{a^*}}\Big)^\frac{2(p-1)}{4-d(p-1)}$, and
$A_\tau$ is chosen  such that $\|u_{\tau}\|_2^2=1$.
Applying the identity \eqref{id:w} and the  exponential decay of $w$ in \eqref{decay:w}, some calculations yield that
\begin{equation*}\label{up:e}
\begin{split}
e(\rho)\leq  E_{\rho}(u_{\tau})
\leq&\frac{1}{2}\frac{d(p-1)}{2(p+1)-d(p-1)}\tau^2
  -\frac{2\big(\frac{\rho}{\sqrt{a^*}}\big)^{p-1}\tau^{\frac{d}{2}(p-1)}}{2(p+1)-d(p-1)}
  +O(1)\\
 =&-(1+o(1))\lambda
 \Big(\frac{\rho}{\sqrt{a^*}}\Big)^\frac{4(p-1)}{4-d(p-1)}
  \,\,\,\text{as}\,\,\,\,\rho\to\infty,\\
\end{split}
\end{equation*}
where $\lambda$ is given in \eqref{def:lam}.
This gives the upper bound of $e(\rho)$ as $\rho\to\infty$.

Next, we shall establish the lower bound estimate of $e(\rho)$ as $\rho\to\infty$
by employing the estimate of $\tilde{e}(\rho)$ given in \eqref{val:tile},
where $\tilde{e}(\rho)$ is a new minimization problem defined by \eqref{def:tile}.
Let $u_\rho$ be a nonnegative minimizer of $e(\rho)$ with $\rho\to\infty$.
Since $\int_{\mathbb{R}^d}V(x)u_\rho^2\dx\geq0$, one can then deduce from \eqref{def:e}, \eqref{def:tile} and \eqref{val:tile} that
\begin{equation*}\label{low:e}
e(\rho)\geq\tilde{e}(\rho)
=-\lambda
 \Big(\frac{\rho}{\sqrt{a^*}}\Big)^\frac{4(p-1)}{4-d(p-1)}\,\,\,\text{as}\,\,\,\rho\to\infty.
\end{equation*}

Combining the upper and lower bound estimates then yields \eqref{val:e}, and  this lemma  is then proved.
\qed

\vskip 0.2truein

\subsection{Blow-up analysis}\label{sec:blowup2}

In this section, we shall complete the proof of Theorem \ref{thm:uk}. 
Let $u_k$ be a nonnegative minimizer of $e(\rho_k)$ with $\rho_k\to\infty$ as $k\to\infty$, and then $u_k$ satisfies \eqref{equ:u} for some suitable $\mu_k$.
We firstly give the following lemma.
\begin{lem}
Suppose $V(x)$ satisfies \eqref{con:V1}.
Let $u_k$ be a nonnegative minimizer of $e(\rho_k)$ with $\rho_k\to\infty$ as $k\to\infty$.
We then have
\begin{equation}\label{lim:e.tile}
 0\leq e(\rho_k)-\tilde{e}(\rho_k)\to 0 \,\,\,\text{as}\,\,\,k\to\infty,
\end{equation}
and
\begin{equation}\label{lim:Vu.0}
\int_{\R^d}V(x)u_k^2\dx\to0\,\,\,\text{as}\,\,\,k\to\infty.
\end{equation}
\end{lem}
As for the proof of this lemma, one can refer to \cite[Lemma 4.2]{M} and we omit it here.
\qed

Define
\begin{equation*}\label{def:bareps}
  \hat{\varepsilon}_k:=\rho_k^{-\frac{2(p-1)}{4-d(p-1)}}
\end{equation*}
and
\begin{equation}\label{def:hatwk}
  \hat{w}_k:=\hat{\varepsilon}_k^{\frac{d}{2}}u_k(\hat{\varepsilon}_k x).
\end{equation}
Some calculations yield that
\begin{equation*}
\int_{\R^d}|\nabla u_k|^2\dx= \hat{\varepsilon}_k^{-2}\int_{\R^d}|\nabla \hat{w}_k|^2\dx
\end{equation*}
and
\begin{equation*}
 \rho_k^{p-1}\int_{\R^d}u_k^{p+1}\dx
=\hat{\varepsilon}_k^{-2}\int_{\R^d}\hat{w}_k^{p+1}\dx.
\end{equation*}
We now give the following lemma on the boundedness of $\int_{\R^d}|\nabla \hat{w}_k|^2\dx$ and $\int_{\R^d}\hat{w}_k^{p+1}\dx$.

\begin{lem}\label{bd:hatwk}
Suppose $V(x)$ satisfies \eqref{con:V1}.
Let $u_k$ be a nonnegative minimizer of $e(\rho_k)$ with $\rho_k\to\infty$ as $k\to\infty$.
Then one has
\begin{equation}\label{bd:u.nabla.p}
 C_1\leq \int_{\R^d}|\nabla \hat{w}_k|^2\dx\leq C_2
 \,\,\,\text{and}\,\,\,
 C'_1\leq \int_{\R^d}\hat{w}_k^{p+1}\dx\leq C'_2,
\end{equation}
where $\hat{w}_k$ is defined by \eqref{def:hatwk}, $C_1$, $C_2$, $C'_1$ and $C'_2$  are positive constants independent of $k$.
\end{lem}
\noindent {\bf Proof.}
It follows from \eqref{def:E} and \eqref{val:e} that
\begin{equation}\label{lim:e.w}
(\sqrt{a^*})^{\frac{4(p-1)}{4-d(p-1)}}\hat{\varepsilon}_k^2e(\rho_k)=\int_{\R^d}|\nabla \hat{w}_k|^2\dx-\frac{1}{p+1}\int_{\R^d}\hat{w}_k^{p+1}\dx
\to-\lambda<0\,\,\,\text{as}\,\,\,k\to\infty,
\end{equation}
where $\lambda$ is given in \eqref{def:lam}.
Using  the Gagliardo-Nirenberg inequality \eqref{ineq:GNw}, one can then derive from \eqref{lim:e.w} that $\int_{\R^d}|\nabla \hat{w}_k|^2\dx\leq C_2$ and $\int_{\R^d}\hat{w}_k^{p+1}\dx\leq C'_2$.
As for the lower bounds,  from \eqref{lim:e.w} one can deduce that $\int_{\R^d}\hat{w}_k^{p+1}\dx\geq C'_1$, and using the Gagliardo-Nirenberg inequality \eqref{ineq:GNw} then yields
$\int_{\R^d}|\nabla \hat{w}_k|^2\dx\geq C\big(\int_{\R^d}\hat{w}_k^{p+1}\dx\big)^\frac{4}{d(p-1)}\geq C_1$.
Hence, we complete the proof of this lemma.
\qed%

\vskip 0.2truein

Motivated by \cite{GWZZ,GZZ,Wang,Zeng1,Zeng2}, we then give the following lemma, which is a weak version of Theorem \ref{thm:uk}.
\begin{lem}\label{lem:wk}
Suppose $V(x)$ satisfies \eqref{con:V1}.
Let $u_k$ be a nonnegative minimizer of $e(\rho_k)$ with $\rho_k\to\infty$ as $k\to\infty$.
We then have follows.
\begin{itemize}
  \item [(i).]
  There exist a sequence $\{y_k\}\subset\R^d$ and positive constants $\iota$ and $R_0$ such that
\begin{equation}\label{subbound:waL2}
 \liminf\limits_{k\to\infty}\int_{B_{R_0}(0)}\hat{w}_k^{p+1}\dx\geq\iota>0.
\end{equation}
  \item [(ii).]
  The sequence $\{y_k\}$ satisfies that, passing to a subsequence if necessary, $\hat{\varepsilon}_ky_k\to z_0$ for some $z_0\in\R^d$ satisfying $V(z_0)=0$.
  \item [(iii).]
  Defined
\begin{equation}\label{def:wk}
  w_k:=\hat{w}_k(x+y_k)=\hat{\varepsilon}_k^{\frac{d}{2}}u_k(\hat{\varepsilon}_k x+\hat{\varepsilon}_k y_k),
\end{equation}
and then passing to a subsequence if necessary, there holds that
 \begin{equation}\label{lim:wk}
  \lim\limits_{k \to \infty}w_k=(\sqrt{a^*})^{-\frac{d(p-1)}{4-d(p-1)}} w\big((\sqrt{a^*})^{-\frac{2(p-1)}{4-d(p-1)}}x+x_0'\big),
\end{equation}
strongly in $H^1(\R^d)$ for some $x_0'\in\R^d$, where $a^*:=\|w\|_2^2$ and $w$  is the unique (up to translations) positive  solution of \eqref{equ:w}.
\end{itemize}
\end{lem}
\noindent {\bf Proof.}
(i). As for \eqref{subbound:waL2}, if it is false, then for any $R > 0$, there exists a subsequence of $\hat{w}_k$  (still denoted by $\hat{w}_k$) such that
$\lim\limits_{k\to\infty}\sup\limits_{y\in\R^d}\int_{B_{R}(y)}\hat{w}_k^{p+1}\dx=0$.
Applying \cite[Lemma 1.1]{Lions2} then yields that $\hat{w}_k\to0$ in $L^{p+1}(\R^d)$ as $k\to\infty$, which however contradicts \eqref{bd:u.nabla.p}.

(ii).
Employing  \eqref{lim:Vu.0}, this conclusion can be obtained by using the proof by contradiction.
Since the proof is similar to that of \cite[Lemma 2.3]{GZZ}, we omit it here.

(iii).
It follows from \eqref{equ:u} and \eqref{def:hatwk} that $w_k$ solves
\begin{equation}\label{equ:wk}
  -\Delta w_k+\hat{\varepsilon}_k^2V(\hat{\varepsilon}_k x+\hat{\varepsilon}_k y_k)w_k=\hat{\varepsilon}_k^2\mu_kw_k+w_k^p.
\end{equation}
Following \eqref{def:E} and \eqref{val:e}, one can deduce that
\begin{equation}\label{lim:mu.wk}
 \hat{\varepsilon}_k^2\mu_k= 2\hat{\varepsilon}_k^2e(\rho_k)-\frac{p-1}{p+1}\int_{\mathbb{R}^d}w_k^{p+1}\dx<0.
\end{equation}
Using  \eqref{bd:u.nabla.p} and \eqref{lim:e.w}, one can then obtain the uniform boundedness of $\{\hat{\varepsilon}_k^2\mu_k\}$ as $k\to\infty$,
which indicates that passing to a subsequence if necessary, $\hat{\varepsilon}_k^2\mu_k\to -\beta$ for some $\beta\in\R^+$ as $k\to\infty$.
From \eqref{bd:u.nabla.p}, one can deduce that $w_k$ is bounded uniformly in $H^1(\R^d)$.
Taking $k\to\infty$,
passing to a subsequence if necessary, one then has  $w_k\rightharpoonup w_0 \geq 0$ in $H^1(\R^d)$ for some $w_0\in H^1(\R^d)$ satisfying
\begin{equation}\label{equ:w0}
  -\Delta w_0+\beta w_0=w_0^p.
\end{equation}
Applying the maximum principle,  one can then conclude from \eqref{subbound:waL2} that $w_0> 0$,
which implies from \eqref{equ:w} that $w_0=\beta^\frac{1}{p-1}w(\beta^\frac{1}{2}x+x_0')$ for some $x_0'\in\R^d$, because
of the uniqueness (up to translations) of positive  solution of \eqref{equ:w}.

Here we claim that
\begin{equation}\label{w0.1}
  \| w_0\|_2^2=1\,\,\,\text{and}\,\,\,\beta=\| w\|_2^{\frac{4(p-1)}{d(p-1)-4}}.
\end{equation}
From the following Pohozaev identity of  \eqref{equ:w0} (cf. \cite[Lemma 8.1.2]{C}),
\[(d-2)\int_{\R^d}|\nabla w_0|^2\dx+d\beta\int_{\R^d}w_0^2\dx=\frac{2d}{p+1}\int_{\R^d}w_0^{p+1}\dx,\]
one can derive that
\begin{equation}\label{id:w0}
\int_{\R^d}|\nabla w_0|^2\dx=\frac{p-1}{p+1}\frac{d}{2}\int_{\R^d}w_0^{p+1}\dx=\frac{d(p-1)}{2(p+1)-d(p-1)}\beta \int_{\R^d}w_0^2\dx.
\end{equation}
Applying the Gagliardo-Nirenberg inequality \eqref{ineq:GNw} and  \eqref{id:w0}, some calculations yield that
\begin{equation*}
\begin{split}
C_{GN}
\leq& \frac{\|\nabla w_0\|_2^{\frac{d}{2}(p-1)}\|w_0\|_2^{p+1-\frac{d}{2}(p-1)}}{\|w_0\|_{p+1}^{p+1}}\\
=&\Big(\frac{p-1}{p+1}\frac{d}{2}\Big)^{\frac{d}{4}(p-1)}
\Big(\frac{2(p+1)\beta}{2(p+1)-d(p-1)}\Big)^{\frac{d}{4}(p-1)-1}\| w_0\|_2^{p-1}\\
\leq&
\frac{C_{GN}}{\| w\|_2^{p-1}}\beta^{\frac{d}{4}(p-1)-1},
\end{split}
\end{equation*}
where $C_{GN}$ is defined in \eqref{ineq:GNw.C}.
This gives that $\beta\geq\| w\|_2^\frac{4(p-1)}{d(p-1)-4}$, and further implies that $\beta=\| w\|_2^\frac{4(p-1)}{d(p-1)-4}$, because $\| w_0\|_2^2= \beta^{\frac{4-d(p-1)}{2(p-1)}}\| w\|^2_2\leq1$.
Moreover, one can deduce that $\| w_0\|_2^2=1$.

Since $\| w_k\|_2=\| w_0\|_2=1$,  passing to a subsequence if necessary, one has $w_k\to w_0 $ strongly in $L^2(\R^d)$ as $k\to\infty$.
Using the interpolation  inequality, one can further derive that $w_k\to w_0 $ strongly in $L^q(\R^d)$ for any $q \in [2,2^*)$ as $k\to\infty$.
Moreover, one can conclude from \eqref{equ:wk} and \eqref{equ:w0}  that (\ref{lim:wk}) holds.
\qed

\vskip 0.2truein

\begin{lem}\label{lem:wk.2}
Suppose $V(x)$ satisfies \eqref{con:V1}.
Let $u_k$ be a nonnegative minimizer of $e(\rho_k)$ with $\rho_k\to\infty$ as $k\to\infty$.
Then $u_k$ admits only one local maximum point $x_k$,
and passing to a subsequence if necessary, there holds that
\begin{equation}\label{lim:barwk.h1}
\bar{u}_k(x):=\sqrt{a^*}\varepsilon_k^\frac{d}{2}u_k(\varepsilon_kx+x_k)\to w(x)\,\,\,\text{strongly in $H^1(\R^d)$ as $k\to\infty$,}
\end{equation}
 where $\varepsilon_k:=\Big(\frac{\rho_k}{\sqrt{a^*}}\Big)^{-\frac{2(p-1)}{4-d(p-1)}}$ is defined in \eqref{def:eps}.
\end{lem}
\noindent {\bf Proof.}
Applying the De-Giorgi-Nash-Moser theory (cf. \cite[Theorem 4.1]{HL}), one can derive from  (\ref{lim:wk}) and \eqref{equ:wk} that
\begin{equation}\label{wkzaza}
   w_k(x)\to 0 \,\, \,\,\text{as}\,\,\, |x| \to \infty \,\, \text{uniformly for large} \,\,k,
\end{equation}
which  indicates that $u_k (x)$ admits at least one global maximum point. Let $x_k $ is a global maximum point of $u_k (x)$ and set $z_k :=\varepsilon_k  y_k \to x_0$ as $k\to\infty$.
Since $z'_k :=\frac{x_k -z_k }{\eps_k }$ is a global maximum point of $w_k (x)$, one can thus derive from \eqref{subbound:waL2} and \eqref{wkzaza} that
\begin{equation}\label{eq3.28}
\Big\{\frac{x_k -z_k }{\eps_k }\Big\} \text{ is  bounded uniformly in } \R^d.
\end{equation}

Define
\begin{equation}\label{def:barwk}
\bar{u}_k(x):=\sqrt{a^*}\varepsilon_k^\frac{d}{2}u_k(\varepsilon_kx+x_k)
\end{equation}
where $\varepsilon_k:=\Big(\frac{\rho_k}{\sqrt{a^*}}\Big)^{-\frac{2(p-1)}{4-d(p-1)}}$ is given in \eqref{def:eps}.
It then follows from \eqref{lim:wk} that, passing to  a subsequence if necessary, $\bar{u}_k(x)\to w(x+y_0')$ for some $y_0'\in\R^d$ strongly in $H^1(\R^d)$ as $k\to\infty$.
Since $V(x)\in C^\alpha(\R^d)$, using the standard elliptic regular theory, we have
\begin{equation}\label{lim:GTP.barwk.C}
\bar{u}_k(x)\to w(x+y_0')\,\,\, \text{in}\,\, C^2_{loc}(\R^d)\, \,\,  \,\text{as}\,\,\, k\to\infty,
\end{equation}
and one can  see \cite[Lemma 3.1]{GZZ} for a detailed proof.
Note that  the origin is a local  maximum point of $\bar{u}_k$ for all $k>0$,
and it follows from \eqref{lim:GTP.barwk.C} that it is also a local  maximum point of $w$.
Since $w(x)$ is radially symmetric about the origin and decreases strictly in $|x|$ (see, e.g., \cite{GNN,LN,W}),
we know that $x = 0$ is the unique local maximum point of $w(x)$,  which thus implies from \eqref{lim:GTP.barwk.C} that $y_0'=0$.
Hence, it follows that
\begin{equation}\label{lim:baruk.H1}
  \bar{u}_k(x)\to w(x)\,\,\,\text{strongly in}\ H^1(\R^d) \, \,\,\text{as}\, \,\, k\to\infty.
\end{equation}

We finally prove the uniqueness of the local maximum points of $u_k$ when $k$ is sufficiently large. Suppose $x_k$ is any local maximum point of $u_k$.
It is easy to know that $\bar{u}_k $ satisfies
\begin{equation}\label{equ:baruk}
   -\Delta \bar{u}_k+\varepsilon_k^2V(\varepsilon_kx+x_k)\bar{u}_k
   =\varepsilon_k^2\mu_k\bar{u}_k+\bar{u}_k^p\,\,\,\text{in}\,\,\,\R^d.
\end{equation}
From this, one can deduce that $\bar{u}_k(x_k)\geq C_0>0$ when $k>0$ is large enough.
This indicates that all local maximum points of $\bar{u}_k $ must stay in a finite ball $B_R(0)$ as $k\to\infty$, where $R>0$ is independent of $k$.
Employing the uniqueness of local maximum points of $w$, one can deduce from \eqref{lim:GTP.barwk.C} that the origin is the unique maximum point of $\bar{u}_k $, $i.e.$,  $u_k$ admits only one local maximum point $x_k$ when $k$ large enough.
\qed

\vskip 0.2truein

\noindent\textbf{Proof of Theorem \ref{thm:uk}. }
As for the exponential decay of $u_k$  in \eqref{decay:u.exp}, one can obtain it by using the comparison principle.
Similar to \eqref{lim:mu.wk}, one can check that $\varepsilon_k^2\mu_k\to-1$ as $k\to\infty$, and then we can derive that, there exists a constant $R>0$ large enough such that
\[-\Delta \bar{u}_k+\frac{1}{2}\bar{u}_k\leq0\,\,\,\text{and}\,\,\,\bar{u}_k\leq Ce^{-\frac{1}{2}R}\,\,\,\text{for}\,\,\,|x|\geq R. \]
Comparing $\bar{u}_k$ with $Ce^{-\frac{1}{2}|x|}$ then yields $\bar{u}_k(x)\le Ce^{-\frac{|x|}{2}}$ when  $|x|\geq R$.
Further more, applying the local elliptic estimate (cf. (3.15) in \cite{GT}) then yields $|\nabla \bar{u}_k|\le Ce^{-\frac{|x|}{4}}$ when  $|x|>R$. We thus give the  proof of  \eqref{decay:u.exp}.
Moreover, by \eqref{decay:u.exp} and \eqref{lim:baruk.H1}, applying the standard elliptic regularity theory then yields \eqref{lim:u},
(see, e.g.,  \cite[Lemma 4.9]{M} for similar arguments).

Finally, we aim at proving \eqref{lim:x}. Suppose $\tilde{u}_{k}$ is a nonnegative minimizer of $\tilde{e}(\rho_{k})$, and then $\tilde{u}_{k}(x-\epsilon_{k} y_{0}-x_{0})$ is also a nonnegative minimizer of $\tilde{e}(\rho_{k})$, where $x_{0}\in Z_{0}$, $y_{0}\in K_{0}$, and $Z_{0}, K_{0}$ are defined by \eqref{def:rho.z0.y0}.
We then derive from \eqref{decay:w}, \eqref{con:V2}, \eqref{con:V3} and \eqref{val:tilu} that
\begin{equation}\label{wwww2}
\begin{split}
  e(\rho_{k})-\tilde{e}(\rho_{k})
  & \leq\int _{\R^d} V(x) \tilde{u}^{2}_{k}(x-\varepsilon_{k} y_{0}-x_{0})\dx\\
  &\leq \frac{1}{a^*}\big(1+o(1)\big)\int_{B_\frac{1}{\sqrt{\varepsilon_{k}}}(0)}
  \frac{V(\varepsilon_{k}x+\varepsilon_{k}y_{0}+x_{0})}{V_{0}(\varepsilon_{k}x+\varepsilon_{k}y_{0})}V_{0}(\varepsilon_{k}x+\varepsilon_{k}y_{0})w^2\dx\\
  &\leq \frac{1}{a^*}\varepsilon^{r}_{k}\big(1+o(1)\big)\int_{\R^d}V_{0}(x+y_{0})w^2\dx
  =\frac{1}{a^*}\varepsilon^{r}_{k}\big(1+o(1)\big)\bar{\lambda}_0,
  \end{split}
\end{equation}
where $\bar{\lambda}_0$ is given by \eqref{def:rho.gamma}.
Suppose $u_{k}$ is a nonnegative minimizer of $e(\rho_{k})$,  and then one can deduce from \eqref{con:V2} and \eqref{decay:u.exp} that
\begin{equation}\label{www1}
\begin{split}
e(\rho_{k})-\tilde{e}(\rho_{k}) &  \geq\int _{\R^d} V(x) u^{2}_{k}\dx =\frac{1}{a^*}\int_{\R^d}V(\varepsilon_kx+x_k)\bar{u}_k^2\dx\\
&=\frac{1}{a^*}\int_{\R^d}\frac{V(\varepsilon_kx+x_k-x_{i}+x_{i})}{V_{i}(\varepsilon_kx+x_k-x_{i})}V_{i}(\varepsilon_kx+x_k-x_{i})\bar{u}_k^2\dx\\
&\geq \frac{1+o(1)}{a^*}\varepsilon^{r_{i}}_{k}\int_{B_{\frac{1}{\sqrt{\varepsilon_{k}}}}(x_{i})}V_{i}\Big(x+\frac{x_k-x_{i}}{\varepsilon_k}\Big)\bar{u}_k^2\dx,\\
\end{split}
\end{equation}
where $x_{i}\in Z$. Comparing with the upper estimate \eqref{wwww2}, one can directly check that $r_i=r$ and $x_i=x_0\in\bar{Z}$, where $r$ and $\bar{Z}$ is given by  \eqref{def:rho.p0}.
Since $V(x)\to\infty$ as $|x|\to\infty$, one can further check that, $\{\frac{x_k-x_{0}}{\varepsilon_k}\}$ is bounded uniformly in $k$.
More precisely, one can also verify that, passing to a subsequence if necessary,
\[
  \frac{x_k-x_{0}}{\varepsilon_k}\to y_0,
\]
which implies that $x_0\in Z_0$ and $Z_0$ is defined in \eqref{def:rho.z0.y0}, $i.e.$, \eqref{lim:x} holds.
Moreover, we also have
\begin{equation}\label{lim:ee}
\lim\limits_{k\to\infty}\frac{e(\rho_k)-\tilde{e}(\rho_k)}{\varepsilon_k^r}
=\frac{1}{a^*}\bar{\lambda}_0,
\end{equation}
where $\bar{\lambda}_0$ is defined by \eqref{def:rho.gamma}.
This gives \eqref{est2:e}, and  the proof of Theorem \ref{thm:uk} is thus completed.
\qed

\section{Local uniqueness of nonnegative minimizers}\label{sec:unique}

In this section, we focus on the proof of local uniqueness of minimizers as $\rho\to\infty$.
Argue by contradiction.  Suppose it is not true, and there exist two different nonnegative minimizers $u_{1k}$ and $u_{2k}$ for $e(\rho_k)$ with $\rho_{k}\to\infty$ as $k\to\infty$. Let $x_{1k}$ and $x_{2k}$ denote the unique local maximum point of $u_{1k}$ and $u_{2k}$, respectively.
Following \eqref{equ:u}, we have
\begin{equation}\label{equ:uik}
   -\Delta u_{ik}+V(x)u_{ik}
   =\mu_{ik}u_{ik}+\rho_k^{p-1}u_{ik}^p\,\,\,\text{in}\,\,\,\R^d,\,\,\,i=1,2.
\end{equation}
Define
\begin{equation}\label{def:baruik}
  \hat{u}_{ik}(x):=\sqrt{a^*}\varepsilon_k^\frac{d}{2}u_{ik}(x)\,\,\,\text{and}\,\,\,
  \bar{u}_{ik}(x):=\hat{u}_{ik}(\varepsilon_kx+x_{1k}),\,\,\, i=1,\ 2,
\end{equation}
where $\varepsilon_k$ is given by \eqref{def:eps}.
Since $\lim\limits_{k\to\infty}\frac{x_{2k}-x_{1k}}{\varepsilon_k}=0$, by Theorem \ref{thm:uk}, one then has $\bar{u}_{ik}\to w(x)$ uniformly in $\R^d$ as $k\to \infty$.
One can check that $\bar{u}_{ik}$ satisfies
\begin{equation}\label{equ:baruik}
   -\Delta \bar{u}_{ik}+\varepsilon_k^2V(\varepsilon_kx+x_{1k})\bar{u}_{ik}
   =\varepsilon_k^2\mu_{ik}\bar{u}_{ik}+\bar{u}_{ik}^p\,\,\,\text{in}\,\,\,\R^d,\,\,\,i=1,2.
\end{equation}
Since $u_{1k}\not\equiv u_{2k}$, define
\begin{equation}\label{def:eta}
  \eta_k:=\frac{ u_{1k}-u_{2k}}{\|u_{1k}- u_{2k}\|_{L^{\infty}(\R^d)}}\,\,\,\text{and}\,\,\,
   \hat{\eta}_k:=\frac{\hat{u}_{1k}-\hat{u}_{2k}}{\|\hat{u}_{1k}- \hat{u}_{2k}\|_{L^{\infty}(\R^d)}},
\end{equation}
and then we have $\eta_k=\hat{\eta}_k$. Further we define
\begin{equation}\label{def:bar.eta}
  \bar{\eta}_k(x):=\hat{\eta}_k(\varepsilon_kx+x_{1k}),
\end{equation}
and thus one can deduce from \eqref{equ:baruik} that  $\bar{\eta}_k$ satisfies
\begin{equation}\label{equ:bar.etak}
   -\Delta \bar{\eta}_k+\varepsilon_k^2V(\varepsilon_kx+x_{1k})\bar{\eta}_k
   =\varepsilon_k^2\mu_{1k}\bar{\eta}_k+\bar{g}_k(x)+\bar{f}_k(x),
\end{equation}
where
\begin{equation}\label{def:gf}
  \bar{g}_k(x):=\varepsilon_k^2\frac{\mu_{1k}-\mu_{2k}}{\|\bar{u}_{1k}- \bar{u}_{2k}\|_{L^{\infty}(\R^d)}}\bar{u}_{2k}\,\,\,\text{and}\,\,\,
  \bar{f}_k(x):=\frac{\bar{u}_{1k}^p-\bar{u}_{2k}^p}{\|\bar{u}_{1k}- \bar{u}_{2k}\|_{L^{\infty}(\R^d)}}.
\end{equation}
Now we give the following lemma on the limit of $\bar{\eta}_k$.

\begin{lem}\label{iem1}
Suppose all the assumptions of Theorem \ref{Thm:uniqueness} hold. Then passing to a subsequence if necessary, $\bar{\eta}_k\to \bar{\eta}_0$ in $C_{loc}(\R^d)$ as $k\to\infty$, where $\bar{\eta}_0$ satisfies
\begin{equation}\label{val:bar.eta}
  \bar{\eta}_0(x)=b_0\big(w+\frac{p-1}{2}x\cdot\nabla w\big)+\sum^d_{i=1}b_i\frac{\partial w}{\partial x_i},
\end{equation}
and $b_0,b_1,...b_d$ are some constants.
\end{lem}
\noindent {\bf Proof.}
Since $\|\bar{\eta}_k\|_\infty\leq1$, the standard elliptic regularity then implies that $\|\bar{\eta}_k\|_{C^{1,\alpha}_{loc}(\R^d)}\leq C$, where $C$ is a constant independent of $k$.
Therefore, passing to a subsequence if necessary, one can deduce that
\begin{equation}\label{lim:etak}
\text{$\bar{\eta}_k\to\bar{\eta}_0$ in $C_{loc}(\R^d)$ as $k\to\infty$ for some function $\bar{\eta}_0\in C_{loc}(\R^d)$.}
\end{equation}
Similar to \eqref{lim:mu.wk}, from \eqref{def:E} and \eqref{def:baruik}, one can derive that
\begin{equation}\label{val:epsmu}
 \varepsilon_k^2\mu_{ik}= 2\varepsilon_k^2e(\rho_k)-\frac{p-1}{a^*(p+1)}\int_{\mathbb{R}^d}|\bar{u}_{ik}|^{p+1}\dx.
\end{equation}
Define
\begin{equation}\label{def:bar.C.p}
\begin{split}
\bar{u}_{1k}^{p+1}-\bar{u}_{2k}^{p+1}
=&\int_{0}^{1}\frac{d}{\dt}\big[t\bar{u}_{1k}+(1-t)\bar{u}_{2k}\big]^{p+1}\dt\\
=&(p+1)(\bar{u}_{1k}-\bar{u}_{2k})\int_{0}^{1}\big[t\bar{u}_{1k}+(1-t)\bar{u}_{2k}\big]^p\dt\\
:=&(p+1)\bar{C}_k^p(x)(\bar{u}_{1k}-\bar{u}_{2k}),
\end{split}
\end{equation}
which implies from \eqref{lim:u} that $\bar{C}_k^p(x)\to w^p(x)$ uniformly in $\R^d$ as $k\to\infty$.
Further one can derive that
\begin{equation}\label{sim:gk}
\begin{split}
\bar{g}_k(x)
=&\varepsilon_k^2\frac{\mu_{1k}-\mu_{2k}}{\|\bar{u}_{1k}- \bar{u}_{2k}\|_{L^{\infty}(\R^d)}}\bar{u}_{2k}
=-\frac{p-1}{a^*(p+1)}
\frac{\int_{\mathbb{R}^d}\big(|\bar{u}_{1k}|^{p+1}-|\bar{u}_{2k}|^{p+1}\big)\dx}{\|\bar{u}_{1k}-\bar{u}_{2k}\|_{L^{\infty}(\R^d)}}
\bar{u}_{2k}\\
=&-\frac{p-1}{a^*}
\int_{\mathbb{R}^d}\bar{C}_k^p(x)\bar{\eta}_k\dx
\bar{u}_{2k},\\
\end{split}
\end{equation}
which implies from \eqref{lim:u} that
\begin{equation}\label{lim:gk}
  \bar{g}_k(x)\to-\frac{p-1}{a^*} \int_{\mathbb{R}^d}w^p\bar{\eta}_0\dx\, w\,\,\,\text{uniformly in}\,\,\,\R^d\,\,\,\text{as}\,\,\,k\to\infty.
\end{equation}
On the other hand, similar to \eqref{def:bar.C.p}, one can also define $\bar{D}_k^{p-1}(x)$ satisfying
\begin{equation}\label{def:bar.C.p-1}
\begin{split}
p\bar{D}_k^{p-1}(x)(\bar{u}_{1k}-\bar{u}_{2k}):=\bar{u}_{1k}^p-\bar{u}_{2k}^p,
\end{split}
\end{equation}
 and then $\bar{D}_k^{p-1}(x)\to w^{p-1}(x)$ uniformly in $\R^d$ as $k\to\infty$.
Further, one has
\begin{equation}\label{sim:fk}
\begin{split}
\bar{f}_k(x)
=\frac{\bar{u}_{1k}^p-\bar{u}_{2k}^p}{\|\bar{u}_{1k}- \bar{u}_{2k}\|_{L^{\infty}(\R^d)}}
=p\bar{D}_k^{p-1}(x)\bar{\eta}_k,
\end{split}
\end{equation}
and
\begin{equation}\label{lim:fk}
  \bar{f}_k(x)\to pw^{p-1}\bar{\eta}_0\,\,\,\text{uniformly in}\,\,\,\R^d\,\,\,\text{as}\,\,\,k\to\infty.
\end{equation}

By the above results, taking $k\to\infty$, it follows from \eqref{equ:bar.etak} that $\bar{\eta}_0$ solves
\begin{equation}\label{equ:bar.eta0}
   -\Delta \bar{\eta}_0+(1-pw^{p-1})\bar{\eta}_0
   =-\frac{p-1}{a^*} \int_{\mathbb{R}^d}w^p\bar{\eta}_0\dx \,w.
\end{equation}
Set $\mathcal{L}:=-\Delta+(1-pw^{p-1})$ and one can check that $\mathcal{L}\big(w+\frac{p-1}{2}x\cdot\nabla w\big)=-(p-1) w$.
Recall from (cf. \cite{Kwong, NT}) that
\begin{equation*}
\ker\mathcal{L}=\Big\{\frac{\partial w}{\partial x_1},\frac{\partial w}{\partial x_2},...\frac{\partial w}{\partial x_d}\Big\},
\end{equation*}
and then one can derive that
\begin{equation}\label{val:eta0}
  \bar{\eta}_0(x)=b_0\Big(w+\frac{p-1}{2}x\cdot\nabla w\Big)+\sum^d_{i=1}b_i\frac{\partial w}{\partial x_i},
\end{equation}
where $b_0$, $b_1$, $b_2$,... and $b_d$ are some constants.
\qed

\vskip 0.2truein

\begin{lem}\label{iem2}
Under the assumptions of Theorem \ref{Thm:uniqueness}, there holds that,
\begin{equation}\label{pid:1}
\begin{split}
\frac{b_0}{2}\frac{p-1}{2}\int_{\R^d}\frac{\partial V_0(x+y_0)}{\partial x_j} \big(x\cdot\nabla w^2\big)
-\sum^d_{i=1}\frac{b_i}{2}\int_{\R^d}\frac{\partial^2 V_0(x+y_0)}{\partial x_i\partial x_j}w^2=0,
\end{split}
\end{equation}
 where $j=1,2,\cdots,d$ and $V_0$ is given by \eqref{def:Q}.
\end{lem}
\noindent {\bf Proof.}
At first, we claim that for any $\bar{x}_0$, there exists a small $\delta>0$ and a constant $C>0$ such that
\begin{equation}\label{est:patB}
\begin{split}
\varepsilon_k^2\int_{\partial B_{\delta}(\bar{x}_0)}|\nabla\hat{\eta}_k|^2\ds
+\varepsilon_k^2\int_{\partial B_{\delta}(\bar{x}_0)}V(x)\hat{\eta}_k^2\ds
+ \int_{\partial B_{\delta}(\bar{x}_0)}\hat{\eta}_k^2\ds
\leq C\varepsilon_k^d,
\end{split}
\end{equation}
where $\hat{\eta}_k$ is given by \eqref{def:eta}.

Following from \eqref{def:eta}, \eqref{def:bar.eta} and \eqref{equ:bar.etak}, one can deduce that $\hat{\eta}_k$ satisfies
\begin{equation}\label{equ:hat.etak}
   -\varepsilon_k^2\Delta \hat{\eta}_k+\varepsilon_k^2V(x)\hat{\eta}_k
   =\varepsilon_k^2\mu_{1k}\hat{\eta}_k+\hat{g}_k(x)+\hat{f}_k(x),
\end{equation}
where
\begin{equation}\label{def:hat.gf}
  \hat{g}_k(x):=\varepsilon_k^2\frac{\mu_{1k}-\mu_{2k}}{\|\hat{u}_{1k}- \hat{u}_{2k}\|_{L^{\infty}(\R^d)}}\hat{u}_{2k}\,\,\,\text{and}\,\,\,
  \hat{f}_k(x):=\frac{\hat{u}_{1k}^p-\hat{u}_{2k}^p}{\|\hat{u}_{1k}- \hat{u}_{2k}\|_{L^{\infty}(\R^d)}}.
\end{equation}
Similar to \eqref{sim:gk} and \eqref{sim:fk}, one has
\begin{equation}\label{sim:hat.gf}
  \hat{g}_k(x)=-\frac{p-1}{a^*}\varepsilon_k^{-d}\int_{\mathbb{R}^d}\hat{C}_k^p(x)\hat{\eta}_k\dx\hat{u}_{2k}\,\,\,\text{and}\,\,\,
  \hat{f}_k(x)=p\hat{D}_k^{p-1}(x)\hat{\eta}_k,
\end{equation}
where $\hat{C}_k^p(\varepsilon_kx+x_{1k}):=\bar{C}_k^p(x)$
 and $\hat{D}_k^{p-1}(\varepsilon_kx+x_{1k}):=\bar{D}_k^{p-1}(x)$.

Multiplying \eqref{equ:hat.etak} by $\hat{\eta}_k$ and integrating over $\R^d$ yield that
\begin{equation*}\label{intequ:hat.eta}
\begin{split}
&\varepsilon_k^2\int_{\R^d}|\nabla\hat{\eta}_k|^2\dx
+\varepsilon_k^2\int_{\R^d}V(x)\hat{\eta}_k^2\dx
-\varepsilon_k^2\mu_{1k} \int_{\R^d}\hat{\eta}_k^2\dx\\
=&p\int_{\R^d}\hat{D}_k^{p-1}(x)\hat{\eta}_k^2\dx
-\frac{p-1}{a^*}\varepsilon_k^{-d}\int_{\mathbb{R}^d}\hat{C}_k^p(x)\hat{\eta}_k\dx\int_{\R^d}\hat{u}_{2k}\hat{\eta}_k\dx\\
=&p\varepsilon_k^d\int_{\R^d}\bar{D}_k^{p-1}(x)\bar{\eta}_k^2\dx
-\frac{p-1}{a^*}\varepsilon_k^d\int_{\mathbb{R}^d}\bar{C}_k^p(x)\bar{\eta}_k\dx\int_{\R^d}\bar{u}_{2k}\bar{\eta}_k\dx\\
=&O(\varepsilon_k^d)\,\,\,\text{as}\,\,\,k\to\infty,
\end{split}
\end{equation*}
where the last equality holds because  $\bar{\eta}_k\to\bar{\eta}_0$, $\bar{u}_{2k}\to w$, $\bar{C}_k^p(x)\to w^p$ and $\bar{D}_k^{p-1}(x)\to w^{p-1}$ uniformly in $\R^d$ as $k\to\infty$.
Applying Lemma 4.5 in \cite{Cao} then yields \eqref{est:patB},
and this completes the proof of this claim.

Following from \eqref{equ:uik}, one can deduce that $\hat{u}_{ik}$ solves
\begin{equation}\label{equ:hat.uik}
   -\varepsilon_k^2\Delta \hat{u}_{ik}+\varepsilon_k^2V(x)\hat{u}_{ik}
   =\varepsilon_k^2\mu_{ik}\hat{u}_{ik}+\hat{u}_{ik}^p\,\,\,\text{in}\,\,\,\R^d,\,\,\,i=1,2.
\end{equation}
Multiplying \eqref{equ:hat.uik} by $\frac{\partial \hat{u}_{ik}}{\partial x_j}$  and integrating over $B_\delta(x_{1k})$, where $i=1,2$, $j=1,2,\cdots,d$ and $\delta$ is given by \eqref{est:patB}, one can obtain the following equality,
\begin{equation}\label{equ:int.uik.1}
\begin{split}
& -\varepsilon_k^2\int_{B_\delta(x_{1k})}\Delta \hat{u}_{ik}\frac{\partial \hat{u}_{ik}}{\partial x_j}
   +\frac{\varepsilon_k^2}{2}\int_{B_\delta(x_{1k})}V(x)\frac{\partial \hat{u}_{ik}^2}{\partial x_j}\\
   =&\frac{\varepsilon_k^2\mu_{ik}}{2}\int_{B_\delta(x_{1k})}\frac{\partial \hat{u}_{ik}^2}{\partial x_j}
   +\frac{1}{p+1}\int_{B_\delta(x_{1k})}\frac{\partial \hat{u}_{ik}^{p+1}}{\partial x_j}.  
\end{split}
\end{equation}

Some calculations yield that
\begin{equation}\label{sim:del.uik.1}
\begin{split}
  & -\varepsilon_k^2\int_{B_\delta(x_{1k})}\Delta \hat{u}_{ik}\frac{\partial \hat{u}_{ik}}{\partial x_j}
=-\varepsilon_k^2\sum_{l=1}^d
\int_{B_\delta(x_{1k})}\frac{\partial^2 \hat{u}_{ik}}{\partial x_l^2}\frac{\partial \hat{u}_{ik}}{\partial x_j}\\
=&-\varepsilon_k^2\sum_{l=1}^d
\Big[
 \int_{\partial B_\delta(x_{1k})}\frac{\partial \hat{u}_{ik}}{\partial x_l}\frac{\partial \hat{u}_{ik}}{\partial x_j}\nu_l\ds
-\int_{B_\delta(x_{1k})}\frac{\partial \hat{u}_{ik}}{\partial x_l}\frac{\partial}{\partial x_j}\frac{\partial \hat{u}_{ik}}{\partial x_l}
\Big]\\
=&-\varepsilon_k^2\sum_{l=1}^d
\Big[
 \int_{\partial B_\delta(x_{1k})}\frac{\partial \hat{u}_{ik}}{\partial x_l}\frac{\partial \hat{u}_{ik}}{\partial x_j}\nu_l\ds
-\frac{1}{2}
 \int_{B_\delta(x_{1k})}\frac{\partial}{\partial x_j}\Big(\frac{\partial\hat{u}_{ik}}{\partial x_l}\Big)^2
\Big]\\
=&-\varepsilon_k^2
\Big[
 \int_{\partial B_\delta(x_{1k})}\frac{\partial \hat{u}_{ik}}{\partial \nu}\frac{\partial \hat{u}_{ik}}{\partial x_j}\ds
-\frac{1}{2}
 \int_{\partial B_\delta(x_{1k})}|\nabla\hat{u}_{ik}|^2\nu_j\ds
\Big],\\
\end{split}
\end{equation}
and
\begin{equation}\label{sim:V.uik.1}
\begin{split}
\frac{\varepsilon_k^2}{2}\int_{B_\delta(x_{1k})}V(x)\frac{\partial \hat{u}_{ik}^2}{\partial x_j}
=&\frac{\varepsilon_k^2}{2}
\Big[
 \int_{\partial B_\delta(x_{1k})}V(x)\hat{u}_{ik}^2\nu_j\ds
-\int_{B_\delta(x_{1k})}\frac{\partial V(x)}{\partial x_j} \hat{u}_{ik}^2
\Big].\\
\end{split}
\end{equation}
It then follows from \eqref{equ:int.uik.1}-\eqref{sim:V.uik.1} that
\begin{equation*}\label{id:Vuik.1}
\begin{split}
&\frac{\varepsilon_k^2}{2}
\int_{B_\delta(x_{1k})}\frac{\partial V(x)}{\partial x_j} \hat{u}_{ik}^2\\
=&\frac{\varepsilon_k^2}{2}
\int_{\partial B_\delta(x_{1k})}V(x)\hat{u}_{ik}^2\nu_j\ds
-\frac{\varepsilon_k^2\mu_{ik}}{2}\int_{\partial B_\delta(x_{1k})}\hat{u}_{ik}^2\nu_j\ds
 -\frac{1}{p+1}\int_{\partial B_\delta(x_{1k})} \hat{u}_{ik}^{p+1}\nu_j\ds
\\
&-\varepsilon_k^2
\Big[
 \int_{\partial B_\delta(x_{1k})}\frac{\partial \hat{u}_{ik}}{\partial \nu}\frac{\partial \hat{u}_{ik}}{\partial x_j}\ds
-\frac{1}{2}
 \int_{\partial B_\delta(x_{1k})}|\nabla\hat{u}_{ik}|^2\nu_j\ds
\Big].\\
\end{split}
\end{equation*}
Further, we have
\begin{equation}\label{id:V.eta.1}
\begin{split}
&\frac{\varepsilon_k^2}{2}
\int_{B_\delta(x_{1k})}\frac{\partial V(x)}{\partial x_j} (\hat{u}_{1k}+\hat{u}_{2k})\hat{\eta}_k\\
=&\frac{\varepsilon_k^2}{2}
\int_{\partial B_\delta(x_{1k})}V(x)(\hat{u}_{1k}+\hat{u}_{2k})\hat{\eta}_k\nu_j\ds
+\frac{\varepsilon_k^2}{2}
 \int_{\partial B_\delta(x_{1k})}(\nabla\hat{u}_{1k}+\nabla\hat{u}_{2k})\nabla\hat{\eta}_k\nu_j\ds\\
&-\frac{\varepsilon_k^2\mu_{1k}}{2}\int_{\partial B_\delta(x_{1k})}(\hat{u}_{1k}+\hat{u}_{2k})\hat{\eta}_k\nu_j\ds
-\frac{1}{2}\int_{\partial B_\delta(x_{1k})}\hat{g}_k\hat{u}_{2k}\nu_j\ds\\
& -\frac{1}{p+1}\int_{\partial B_\delta(x_{1k})}\frac{\hat{u}_{1k}^{p+1}-\hat{u}_{2k}^{p+1}}{\|\hat{u}_{1k}- \hat{u}_{2k}\|_{L^{\infty}(\R^d)}}\nu_j\ds\\
&-\varepsilon_k^2
\Big[
 \int_{\partial B_\delta(x_{1k})}\frac{\partial \hat{u}_{1k}}{\partial \nu}\frac{\partial \hat{\eta}_k}{\partial x_j}\ds
+ \int_{\partial B_\delta(x_{1k})}\frac{\partial \hat{\eta}_k}{\partial \nu}\frac{\partial \hat{u}_{2k}}{\partial x_j}\ds
\Big],\\
\end{split}
\end{equation}
where $\hat{g}_k$ is defined in \eqref{def:hat.gf}.

Using the H\"{o}lder inequality, one can derive from \eqref{est:patB} and \eqref{decay:u.exp} that
\begin{equation}\label{est:1.1}
\begin{split}
& \varepsilon_k^2\Big|\int_{\partial B_\delta(x_{1k})}\frac{\partial \hat{u}_{1k}}{\partial \nu}\frac{\partial \hat{\eta}_k}{\partial x_j}\ds
+ \int_{\partial B_\delta(x_{1k})}\frac{\partial \hat{\eta}_k}{\partial \nu}\frac{\partial \hat{u}_{2k}}{\partial x_j}\ds\Big|\\
\leq &\varepsilon_k^2 \Big(\int_{\partial B_\delta(x_{1k})}|\nabla\hat{\eta}_k|^2\Big)^\frac{1}{2}
\Big[\Big(\int_{\partial B_\delta(x_{1k})}|\nabla\hat{u}_{1k}|^2\Big)^\frac{1}{2}
+\Big(\int_{\partial B_\delta(x_{1k})}|\nabla\hat{u}_{2k}|^2\Big)^\frac{1}{2}
\Big]\\
\leq & C\varepsilon_k\varepsilon_k^{\frac{d}{2}} \varepsilon_k^{\frac{d-3}{2}}e^{-\frac{c\delta}{\varepsilon_k}}
=C\varepsilon_k^{d-\frac{1}{2}}e^{-\frac{c\delta}{\varepsilon_k}}
=o(1)e^{-\frac{c\delta}{\varepsilon_k}},
\\
\end{split}
\end{equation}
where $C$ is a suitable positive constant. Similarly, we also have 
\begin{equation}\label{est:1.2}
\begin{split}
\Big|\frac{\varepsilon_k^2}{2}
 \int_{\partial B_\delta(x_{1k})}(\nabla\hat{u}_{1k}+\nabla\hat{u}_{2k})\nabla\hat{\eta}_k\nu_j\ds\Big|
=o(1)e^{-\frac{c\delta}{\varepsilon_k}},
\end{split}
\end{equation}
\begin{equation}\label{est:1.3}
\begin{split}
&\frac{\varepsilon_k^2}{2}
\Big|\int_{\partial B_\delta(x_{1k})}V(x)(\hat{u}_{1k}+\hat{u}_{2k})\hat{\eta}_k\nu_j\ds\Big|\\
\leq&\frac{\varepsilon_k^2}{2}
\Big(\int_{\partial B_\delta(x_{1k})}V(x)\hat{\eta}_k^2\ds\Big)^\frac{1}{2}
\Big(\int_{\partial B_\delta(x_{1k})}V(x)(\hat{u}_{1k}+\hat{u}_{2k})^2\ds\Big)^\frac{1}{2}
\leq o(1)e^{-\frac{c\delta}{\varepsilon_k}},\\
\end{split}
\end{equation}
\begin{equation}\label{est:1.4}
\begin{split}
&\Big|\frac{\varepsilon_k^2\mu_{ik}}{2}\int_{\partial B_\delta(x_{1k})}(\hat{u}_{1k}+\hat{u}_{2k})\hat{\eta}_k\nu_j\ds\Big|
+\frac{1}{2}\Big|\int_{\partial B_\delta(x_{1k})}\hat{g}_k\hat{u}_{2k}\nu_j\ds\Big|\\
\leq&C
\Big(\int_{\partial B_\delta(x_{1k})}(\hat{u}_{1k}+\hat{u}_{2k})^2\ds\Big)^\frac{1}{2}
\Big(\int_{\partial B_\delta(x_{1k})}\hat{\eta}_k^2\ds\Big)^\frac{1}{2}\\
&+C\int_{\mathbb{R}^d}|\bar{C}_k^p(x)|\bar{\eta}_k\dx\int_{\partial B_\delta(x_{1k})}\hat{u}_{2k}^2\ds\\
=&o(1)e^{-\frac{c\delta}{\varepsilon_k}},\\
\end{split}
\end{equation}
and
\begin{equation}\label{est:1.5}
\Big|\frac{1}{p+1}\int_{\partial B_\delta(x_{1k})}\frac{\hat{u}_{1k}^{p+1}-\hat{u}_{2k}^{p+1}}{\|\hat{u}_{1k}- \hat{u}_{2k}\|_{L^{\infty}(\R^d)}}\nu_j\Big|
\leq C\int_{\partial B_\delta(x_{1k})}|\hat{C}_k^p(x)|\hat{\eta}_k
=o(1)e^{-\frac{c\delta}{\varepsilon_k}}.
\end{equation}

On the other hand, let $x_0$ be the unique point of $Z_0$, where $Z_0$ satisfies \eqref{def:rho.z0.y0} and \eqref{con:V4}.
Employing  \eqref{est:1.1}-\eqref{est:1.5}, and
applying \eqref{con:V2}, \eqref{con:V5} and \eqref{lim:etak},  one can derive from \eqref{id:V.eta.1} that
\begin{equation}\label{est:1.6}
\begin{split}
o(1)e^{-\frac{c\delta}{\varepsilon_k}}
=&\frac{\varepsilon_k^2}{2}
\int_{B_\delta(x_{1k})}\frac{\partial V(x)}{\partial x_j} (\hat{u}_{1k}+\hat{u}_{2k})\hat{\eta}_k\\
=&\varepsilon_k^d\frac{\varepsilon_k^2}{2}
\int_{B_\frac{\delta}{\varepsilon_k}(0)}\frac{\partial V(\varepsilon_kx+x_{1k}-x_{0}+x_{0})}{\partial \varepsilon_kx_j} (\bar{u}_{1k}+\bar{u}_{2k})\bar{\eta}_k\\
=&\frac{\varepsilon_k^{d+r+1}}{2}
\int_{B_\frac{\delta}{\varepsilon_k}(0)}\frac{\partial V_{0}(x+\frac{x_{1k}-x_{0}}{\varepsilon_k})}{\partial x_j} (\bar{u}_{1k}+\bar{u}_{2k})\bar{\eta}_k\\
&\quad+\varepsilon_k^{d+2}\int_{B_\frac{\delta}{\varepsilon_k}(0)} W_{j}(\varepsilon_kx+x_{1k}-x_{0}) (\bar{u}_{1k}+\bar{u}_{2k})\bar{\eta}_k\\
=&(1+o(1))\varepsilon_k^{d+r+1}
\int_{\R^d}\frac{\partial V_{0}(x+y_0)}{\partial x_j} w\bar{\eta}_0.
\\
\end{split}
\end{equation}
Furthermore, one can deduce from \eqref{val:eta0} and \eqref{est:1.6} that
\begin{equation*}\label{sim:V.eta0.1}
\begin{split}
0
=&\int_{\R^d}\frac{\partial V_0(x+y_0)}{\partial x_j} w\bar{\eta}_0\\
=&\int_{\R^d}\frac{\partial V_0(x+y_0)}{\partial x_j} w
\Big[b_0\big(w+\frac{p-1}{2}x\cdot\nabla w\big)+\sum^d_{i=1}b_i\frac{\partial w}{\partial x_i}\Big]\\
=&b_0\int_{\R^d}\frac{\partial V_0(x+y_0)}{\partial x_j} w^2
+\frac{b_0}{2}\frac{p-1}{2}\int_{\R^d}\frac{\partial V_0(x+y_0)}{\partial x_j} x\cdot\nabla w^2\\
&+\sum^d_{i=1}\frac{b_i}{2}\int_{\R^d}\frac{\partial V_0(x+y_0)}{\partial x_j}\frac{\partial w^2}{\partial x_i}\\
=&\frac{b_0}{2}\frac{p-1}{2}\int_{\R^d}\frac{\partial V_0(x+y_0)}{\partial x_j}\big( x\cdot\nabla w^2\big)
-\sum^d_{i=1}\frac{b_i}{2}\int_{\R^d}\frac{\partial^2 V_0(x+y_0)}{\partial x_i\partial x_j}w^2,\\
\end{split}
\end{equation*}
which gives \eqref{pid:1}.
\qed
\vskip 0.2truein

In the following, we shall follow the above two lemma to complete the proof of Theorem \ref{Thm:uniqueness}.

\noindent\textbf{Proof of Theorem \ref{Thm:uniqueness}.}
At first, we claim that the coefficient $b_0$ given in \eqref{val:eta0} satisfies
\begin{equation}\label{val:b0.0}
b_0\equiv0.
\end{equation}
Multiplying \eqref{equ:hat.uik} by $(x-x_{1k})\cdot\nabla  \hat{u}_{ik}$ and integrating over $B_\delta(x_{1k})$, where $i=1,2$ and $\delta$ is given in \eqref{est:patB}, one has
\begin{equation}\label{equ:id.2}
\begin{split}
   &-\varepsilon_k^2\int_{B_\delta(x_{1k})}\Delta \hat{u}_{ik}\big[(x-x_{1k})\cdot\nabla  \hat{u}_{ik}\big]
   +\frac{\varepsilon_k^2}{2}\int_{B_\delta(x_{1k})}V(x)\big[(x-x_{1k})\cdot\nabla  \hat{u}_{ik}^2\big]\\
   =&\frac{\varepsilon_k^2\mu_{ik}}{2}\int_{B_\delta(x_{1k})}\big[(x-x_{1k})\cdot\nabla  \hat{u}_{ik}^2\big]
   +\frac{1}{p+1}\int_{B_\delta(x_{1k})}\big[(x-x_{1k})\cdot\nabla  \hat{u}_{ik}^{p+1}\big]. \\
\end{split}
\end{equation}

Some calculations yield that
\begin{equation}\label{sim.del.2}
\begin{split}
&-\varepsilon_k^2\int_{B_\delta(x_{1k})}\Delta \hat{u}_{ik}[(x-x_{1k})\cdot\nabla \hat{u}_{ik}]\\
=&-\varepsilon_k^2\int_{\partial B_\delta(x_{1k})} \frac{\partial\hat{u}_{ik}}{\partial \nu}[(x-x_{1k})\cdot\nabla \hat{u}_{ik}]
  +\varepsilon_k^2\int_{B_\delta(x_{1k})}\nabla \hat{u}_{ik}\cdot\nabla[(x-x_{1k})\cdot\nabla \hat{u}_{ik}]\\
=&-\varepsilon_k^2\int_{\partial B_\delta(x_{1k})} \frac{\partial\hat{u}_{ik}}{\partial \nu}[(x-x_{1k})\cdot\nabla \hat{u}_{ik}]
  +\frac{\varepsilon_k^2}{2}\int_{\partial B_\delta(x_{1k})}[(x-x_{1k})\cdot\nu]|\nabla \hat{u}_{ik}|^2\\
&+\frac{2-d}{4}\varepsilon_k^2\int_{\partial B_\delta(x_{1k})}\nabla\hat{u}_{ik}^2\cdot\nu\\
&-\frac{2-d}{2}
\Big[\varepsilon_k^2\int_{B_\delta(x_{1k})}V(x)\hat{u}_{ik}^2
-\varepsilon_k^2\mu_{ik}\int_{B_\delta(x_{1k})}\hat{u}_{ik}^2
-\int_{B_\delta(x_{1k})}\hat{u}_{ik}^{p+1}\Big],\\
\end{split}
\end{equation}
where the last $``="$ holds due to that
\begin{equation*}
\begin{split}
&\varepsilon_k^2\int_{B_\delta(x_{1k})}\nabla \hat{u}_{ik}\cdot\nabla[(x-x_{1k})\cdot\nabla \hat{u}_{ik}]\\
=&\varepsilon_k^2\sum_{j=1}^d\int_{B_\delta(x_{1k})} \frac{\partial\hat{u}_{ik}}{\partial x_j}
\Big[\frac{\partial\hat{u}_{ik}}{\partial x_j}+(x-x_{1k})\cdot\nabla \frac{\partial\hat{u}_{ik}}{\partial x_j}\Big]\\
=&\varepsilon_k^2\sum_{j=1}^d\int_{B_\delta(x_{1k})} \Big(\frac{\partial\hat{u}_{ik}}{\partial x_j}\Big)^2
+\frac{\varepsilon_k^2}{2}\sum_{j=1}^d\int_{B_\delta(x_{1k})} (x-x_{1k})\cdot\nabla \Big(\frac{\partial\hat{u}_{ik}}{\partial x_j}\Big)^2\\
=&\varepsilon_k^2\int_{B_\delta(x_{1k})}|\nabla\hat{u}_{ik}|^2
+\frac{\varepsilon_k^2}{2}\int_{B_\delta(x_{1k})} (x-x_{1k})\cdot\nabla|\nabla\hat{u}_{ik}|^2\\
=&\frac{\varepsilon_k^2}{2}\int_{\partial B_\delta(x_{1k})} |\nabla\hat{u}_{ik}|^2(x-x_{1k})\cdot\nu
+\frac{2-d}{2}\varepsilon_k^2\int_{B_\delta(x_{1k})}|\nabla\hat{u}_{ik}|^2,\\
\end{split}
\end{equation*}
and
\begin{equation*}
\begin{split}
&\frac{2-d}{2}\varepsilon_k^2\int_{B_\delta(x_{1k})}|\nabla\hat{u}_{ik}|^2\\
=&\frac{2-d}{4}\varepsilon_k^2\int_{\partial B_\delta(x_{1k})}\nabla\hat{u}_{ik}^2\cdot\nu
-\frac{2-d}{2}\varepsilon_k^2\int_{B_\delta(x_{1k})}\hat{u}_{ik}\Delta\hat{u}_{ik}\\
=&\frac{2-d}{4}\varepsilon_k^2\int_{\partial B_\delta(x_{1k})}\nabla\hat{u}_{ik}^2\cdot\nu\\
&-\frac{2-d}{2}
\Big[\varepsilon_k^2\int_{B_\delta(x_{1k})}V(x)\hat{u}_{ik}^2
-\varepsilon_k^2\mu_{ik}\int_{B_\delta(x_{1k})}\hat{u}_{ik}^2
-\int_{B_\delta(x_{1k})}\hat{u}_{ik}^{p+1}\Big].\\
\end{split}
\end{equation*}
Moreover, one can also deduce that
\begin{equation*}\label{sim.mu.2}
\begin{split}
&\int_{B_\delta(x_{1k})}[(x-x_{1k})\cdot\nabla  \hat{u}_{ik}^2]
=\int_{\partial B_\delta(x_{1k})}\hat{u}_{ik}^2[(x-x_{1k})\cdot\nu]
-d\int_{B_\delta(x_{1k})}\hat{u}_{ik}^2,\\
\end{split}
\end{equation*}
\begin{equation*}\label{sim.up.2}
\begin{split}
&\int_{B_\delta(x_{1k})}[(x-x_{1k})\cdot\nabla  \hat{u}_{ik}^{p+1}]
=\int_{\partial B_\delta(x_{1k})}\hat{u}_{ik}^{p+1}[(x-x_{1k})\cdot\nu]
-d\int_{B_\delta(x_{1k})}\hat{u}_{ik}^{p+1},\\
\end{split}
\end{equation*}
and
\begin{equation*}\label{sim.Vu.2}
\begin{split}
&\int_{B_\delta(x_{1k})}V(x)[(x-x_{1k})\cdot\nabla  \hat{u}_{ik}^2]\\
=&\int_{\partial B_\delta(x_{1k})}V(x)\hat{u}_{ik}^2[(x-x_{1k})\cdot\nu]
-\int_{B_\delta(x_{1k})}\big[\nabla V(x)\cdot(x-x_{1k})\hat{u}_{ik}^2
+d V(x)\hat{u}_{ik}^2\big].\\
\end{split}
\end{equation*}

Substituting the above results into \eqref{equ:id.2} yields that
\begin{equation*}
\begin{split}
&-\varepsilon_k^2\int_{\partial B_\delta(x_{1k})} \frac{\partial\hat{u}_{ik}}{\partial \nu}[(x-x_{1k})\cdot\nabla \hat{u}_{ik}]
  +\frac{\varepsilon_k^2}{2}\int_{\partial B_\delta(x_{1k})}[(x-x_{1k})\cdot\nu]|\nabla \hat{u}_{ik}|^2\\
&+\frac{2-d}{4}\varepsilon_k^2\int_{\partial B_\delta(x_{1k})}\nabla\hat{u}_{ik}^2\cdot\nu\\
&-\frac{2-d}{2}
\Big[\varepsilon_k^2\int_{B_\delta(x_{1k})}V(x)\hat{u}_{ik}^2
-\varepsilon_k^2\mu_{ik}\int_{B_\delta(x_{1k})}\hat{u}_{ik}^2
-\int_{B_\delta(x_{1k})}\hat{u}_{ik}^{p+1}\Big]\\
  &-\frac{d\varepsilon_k^2}{2}\int_{B_\delta(x_{1k})}V(x)\hat{u}_{ik}^2-\frac{\varepsilon^{2}_{k}}{2}\int_{B_\delta(x_{1k})}\nabla V(x)\cdot(x-x_{1k})\hat{u}_{ik}^2
\\
&+\frac{\varepsilon_k^2}{2}\int_{\partial B_\delta(x_{1k})}V(x)\hat{u}_{ik}^2[(x-x_{1k})\cdot\nu]\\
=&\frac{\varepsilon_k^2\mu_{ik}}{2}\int_{\partial B_\delta(x_{1k})}\hat{u}_{ik}^2[(x-x_{1k})\cdot\nu]
-\frac{d\varepsilon_k^2\mu_{ik}}{2}\int_{B_\delta(x_{1k})}\hat{u}_{ik}^2\\
&+\frac{1}{p+1}\int_{\partial B_\delta(x_{1k})}\hat{u}_{ik}^{p+1}[(x-x_{1k})\cdot\nu]
-\frac{d}{p+1}\int_{B_\delta(x_{1k})}\hat{u}_{ik}^{p+1}.\\
\end{split}
\end{equation*}
Following from \eqref{val:epsmu}, one has
\begin{equation*}
\varepsilon_k^2\mu_{ik}\int_{\R^d}\hat{u}_{ik}^2= 2a^*\varepsilon_k^d\varepsilon_k^2e(\rho_k)-\frac{p-1}{p+1}\int_{\mathbb{R}^d}|\hat{u}_{ik}|^{p+1}\dx,
\end{equation*}
and it then follows that
\begin{equation}\label{id:Vu.2}
\begin{split}
&-\varepsilon_k^2\int_{\partial B_\delta(x_{1k})} \frac{\partial\hat{u}_{ik}}{\partial \nu}[(x-x_{1k})\cdot\nabla \hat{u}_{ik}]
  +\frac{\varepsilon_k^2}{2}\int_{\partial B_\delta(x_{1k})}[(x-x_{1k})\cdot\nu]|\nabla \hat{u}_{ik}|^2\\
&+\frac{2-d}{4}\varepsilon_k^2\int_{\partial B_\delta(x_{1k})}\nabla\hat{u}_{ik}^2\cdot\nu\\
&+\frac{\varepsilon_k^2}{2}\int_{\partial B_\delta(x_{1k})}V(x)\hat{u}_{ik}^2[(x-x_{1k})\cdot\nu]
-\frac{\varepsilon_k^2}{2}\int_{B_\delta(x_{1k})}[\nabla V(x)\cdot (x-x_{1k})]\hat{u}_{ik}^2\\
&-\frac{\varepsilon_k^2\mu_{ik}}{2}\int_{\partial B_\delta(x_{1k})}\hat{u}_{ik}^2[(x-x_{1k})\cdot\nu]
-\frac{1}{p+1}\int_{\partial B_\delta(x_{1k})}\hat{u}_{ik}^{p+1}[(x-x_{1k})\cdot\nu]\\
=&\varepsilon_k^2       \int_{B_\delta(x_{1k})}V(x)\hat{u}_{ik}^2
-\varepsilon_k^2\mu_{ik}\int_{B_\delta(x_{1k})}\hat{u}_{ik}^2
-\frac{2(p+1)-d(p-1)}{2(p+1)}        \int_{B_\delta(x_{1k})}\hat{u}_{ik}^{p+1}\\
=&\varepsilon_k^2       \int_{B_\delta(x_{1k})}V(x)\hat{u}_{ik}^2
+\frac{d(p-1)-4}{2(p+1)}        \int_{\R^d}\hat{u}_{ik}^{p+1}
-2a^*\varepsilon_k^d\varepsilon_k^2e(\rho_k)\\
&+\varepsilon_k^2\mu_{ik}\int_{\R^d\setminus B_\delta(x_{1k})}\hat{u}_{ik}^2
+\frac{2(p+1)-d(p-1)}{2(p+1)}        \int_{\R^d\setminus B_\delta(x_{1k})}\hat{u}_{ik}^{p+1}.\\
\end{split}
\end{equation}

By \eqref{def:eta}, one can deduce from \eqref{id:Vu.2} that
\begin{equation}\label{id:Veta.2}
\begin{split}
\frac{d(p-1)-4}{2(p+1)}\int_{\R^d}\frac{\hat{u}_{1k}^{p+1}-\hat{u}_{2k}^{p+1}}{\|\hat{u}_{1k}- \hat{u}_{2k}\|_{L^{\infty}(\R^d)}}
=T_1+T_2+T_3+T_4+T_5,
\end{split}
\end{equation}
where
\begin{equation*}\label{def:id2.T1}
\begin{split}
T_1:=&-\varepsilon_k^2\int_{\partial B_\delta(x_{1k})} \frac{\partial\hat{u}_{1k}}{\partial \nu}[(x-x_{1k})\cdot\nabla \hat{\eta}_k]
-\varepsilon_k^2\int_{\partial B_\delta(x_{1k})} \frac{\partial\hat{\eta}_k}{\partial \nu}[(x-x_{1k})\cdot\nabla \hat{u}_{2k}]\\
&  +\frac{\varepsilon_k^2}{2}\int_{\partial B_\delta(x_{1k})}[(x-x_{1k})\cdot\nu]
\nabla \hat{\eta}_k\cdot\nabla(\hat{u}_{1k}+\hat{u}_{2k})\\
&+\frac{2-d}{4}\varepsilon_k^2
\Big[\int_{\partial B_\delta(x_{1k})}\hat{\eta}_k\big(\nabla(\hat{u}_{1k}+\hat{u}_{2k})\cdot\nu\big)
+\int_{\partial B_\delta(x_{1k})}(\hat{u}_{1k}+\hat{u}_{2k})\big(\nabla\hat{\eta}_k\cdot\nu\big)\Big],\\
\end{split}
\end{equation*}
\begin{equation*}\label{def:id2.T2}
\begin{split}
T_2:=&\frac{\varepsilon_k^2}{2}\int_{\partial B_\delta(x_{1k})}V(x)[(x-x_{1k})\cdot\nu](\hat{u}_{1k}+\hat{u}_{2k})\hat{\eta}_k,\\
\end{split}
\end{equation*}
\begin{equation*}\label{def:id2.T3}
\begin{split}
T_3:=
&-\frac{\varepsilon_k^2\mu_{1k}}{2}\int_{\partial B_\delta(x_{1k})}[(x-x_{1k})\cdot\nu](\hat{u}_{1k}+\hat{u}_{2k})\hat{\eta}_k
-\frac{1}{2}\int_{\partial B_\delta(x_{1k})}[(x-x_{1k})\cdot\nu]\hat{g}_k\hat{u}_{2k}\\
&-\varepsilon_k^2\mu_{1k}\int_{\R^d\setminus B_\delta(x_{1k})}(\hat{u}_{1k}+\hat{u}_{2k})\hat{\eta}_k
-\int_{\R^d\setminus B_\delta(x_{1k})}\hat{g}_k\hat{u}_{2k},\\
\end{split}
\end{equation*}
\begin{equation*}\label{def:id2.T4}
\begin{split}
T_4:=
&-\frac{1}{p+1}\int_{\partial B_\delta(x_{1k})}[(x-x_{1k})\cdot\nu]
\frac{\hat{u}_{1k}^{p+1}-\hat{u}_{2k}^{p+1}}{\|\hat{u}_{1k}- \hat{u}_{2k}\|_{L^{\infty}(\R^d)}}\\
&-\frac{2(p+1)-d(p-1)}{2(p+1)}  \int_{\R^d\setminus  B_\delta(x_{1k})}\frac{\hat{u}_{1k}^{p+1}-\hat{u}_{2k}^{p+1}}{\|\hat{u}_{1k}- \hat{u}_{2k}\|_{L^{\infty}(\R^d)}},\\
\end{split}
\end{equation*}
and
\begin{equation*}\label{def:id2.T5}
T_5:=-\varepsilon_k^2\int_{B_\delta(x_{1k})}V(x)(\hat{u}_{1k}+\hat{u}_{2k})\hat{\eta}_k
-\frac{\varepsilon_k^2}{2}\int_{B_\delta(x_{1k})}[\nabla V(x)\cdot (x-x_{1k})](\hat{u}_{1k}+\hat{u}_{2k})\hat{\eta}_k.
\end{equation*}

Similar to the estimates \eqref{est:1.1}-\eqref{est:1.5}, 
one can deduce that
\begin{equation}\label{est:id2.T134}
  |T_1|,|T_2|, |T_3|, |T_4| =o(1)e^{-\frac{c\delta}{\varepsilon_k}}.
\end{equation}
As for $T_5$, the estimate \eqref{est:1.6} gives that
\begin{equation*}
  \frac{\varepsilon_k^2}{2}\Big|\int_{B_\delta(x_{1k})}[\nabla V(x)\cdot x_{1k}](\hat{u}_{1k}+\hat{u}_{2k})\hat{\eta}_k\Big|
  =o(1)e^{-\frac{c\delta}{\varepsilon_k}},
\end{equation*}
and by \eqref{con:V2}, one has
\begin{equation*}
\begin{split}
&\varepsilon_k^2\int_{B_\delta(x_{1k})}V(x)(\hat{u}_{1k}+\hat{u}_{2k})\hat{\eta}_k\\
=&\varepsilon_k^{2+d}\int_{B_\frac{\delta}{\varepsilon_k}(0)}
\frac{V(\varepsilon_kx+x_{1k}-x_{0}+x_{0})}{V_0(\varepsilon_kx+x_{1k}-x_{0})}V_0(\varepsilon_kx+x_{1k}-x_{0})(\bar{u}_{1k}+\bar{u}_{2k})\bar{\eta}_k\\
=&(1+o(1))\varepsilon_k^{2+d+r}\int_{B_\frac{\delta}{\varepsilon_k}(0)}
V_0\Big(x+\frac{x_{1k}-x_{0}}{\varepsilon_k}\Big)(\bar{u}_{1k}+\bar{u}_{2k})\bar{\eta}_k\\
=&(2+o(1))\varepsilon_k^{2+r+d}\int_{\R^d}V_0(x+y_0)w\bar{\eta}_0,\\
\end{split}
\end{equation*}
where $V_0$ is given by \eqref{def:Q}  and $x_0\in Z_0$ with $Z_0$ defined by \eqref{def:rho.z0.y0}.
Moreover, since $\nabla V_{0}(x)\cdot x=r V_{0}(x)$, one can derive from \eqref{con:V5} and \eqref{est:1.6} that
\begin{equation*}\label{qqq2}
\begin{split}
  &\frac{\varepsilon_k^2}{2}\int_{B_\delta(x_{1k})}\big[\nabla V(x)\cdot x\big](\hat{u}_{1k}+\hat{u}_{2k})\hat{\eta}_k\\
  =&\frac{\varepsilon_k^2}{2}\int_{B_\delta(x_{1k})}\big[\nabla V(x-x_{0}+x_{0})\cdot (x-x_{0}+x_{0})\big](\hat{u}_{1k}+\hat{u}_{2k})\hat{\eta}_k\\
  =&\frac{\varepsilon_k^2}{2}\int_{B_\delta(x_{1k})}\big[\nabla V_{0}(x-x_{0})+W(x-x_{0})\big]\cdot (x-x_{0})(\hat{u}_{1k}+\hat{u}_{2k})\hat{\eta}_k\\
  &+\frac{\varepsilon_k^2}{2}\int_{B_\delta(x_{1k})}\big[\nabla V_0(x)\cdot x_{0}\big](\hat{u}_{1k}+\hat{u}_{2k})\hat{\eta}_k\\
  =&\frac{\varepsilon_k^2}{2}\int_{B_\delta(x_{1k})}\big[rV_{0}(x-x_{0})+W(x-x_{0})\big(x-x_{0}\big)\big](\hat{u}_{1k}+\hat{u}_{2k})\hat{\eta}_k
  +o(1)e^{-\frac{c\delta}{\varepsilon_k}}\\
  =&\frac{r(1+o(1))}{2}\varepsilon_k^{d+r+1}\int_{B_\frac{\delta}{\varepsilon_k}(0)}
   V_{0}\Big(x+\frac{x_{1k}-x_{0}}{\varepsilon_{k}}\Big)
  (\bar{u}_{1k} +\bar{u}_{2k})\bar{\eta}_k +o(1)e^{-\frac{c\delta}{\varepsilon_k}}\\
  =&(1+o(1))r\varepsilon_k^{2+r+d}\int_{\R^d}V_0(x+y_0)w\bar{\eta}_0,
  \end{split}
\end{equation*}
where $W:=(W_1,W_2,\cdots,W_d)$.
It then follows from these estimates that
\begin{equation}\label{est:id2.T5}
 T_5=O(1)\varepsilon_k^{2+r+d}.
\end{equation}

As for the left hand of \eqref{id:Veta.2}, one can deduce from \eqref{est:id2.T134} and \eqref{est:id2.T5} that
\begin{equation}\label{decay:id.2}
\begin{split}
O(1)\varepsilon_k^{2+r+d}
=&\frac{d(p-1)-4}{2(p+1)}\int_{\R^d}\frac{\hat{u}_{1k}^{p+1}-\hat{u}_{2k}^{p+1}}{\|\hat{u}_{1k}- \hat{u}_{2k}\|_{L^{\infty}(\R^d)}}\\
=&\frac{d(p-1)-4}{2}\varepsilon_k^d\int_{\R^d}\bar{C}_k^p(x)\bar{\eta}_k\\
=&\frac{d(p-1)-4}{2}(1+o(1))\varepsilon_k^d\int_{\R^d}w^p\bar{\eta}_0,
\\
\end{split}
\end{equation}
where $\bar{C}_k^p(x)$ is defined in \eqref{def:bar.C.p} and the last $``="$ holds because $\bar{C}_k^p(x) \to w^p$  uniformly in $\R^d$ as $k\to\infty$. 
Using \eqref{val:eta0}, then \eqref{decay:id.2} gives that
\begin{equation*}\label{sim:V.eta0.2}
\begin{split}
0
=&\int_{\R^d}w^p\bar{\eta}_0\\
=&\int_{\R^d}w^p\Big[b_0\big(w+\frac{p-1}{2}x\cdot\nabla w\big)+\sum^d_{i=1}b_i\frac{\partial w}{\partial x_i}\Big]\\
=&b_0\int_{\R^d}w^{p+1}+\frac{b_0}{2}\frac{p-1}{p+1}\int_{\R^d}x\cdot\nabla w^{p+1}
+\sum^d_{i=1}\frac{b_i}{p+1}\int_{\R^d}\frac{\partial w^{p+1}}{\partial x_i}\\
=&b_0\int_{\R^d}w^{p+1}-\frac{db_0}{2}\frac{p-1}{p+1}\int_{\R^d} w^{p+1}\\
=&\Big[1-\frac{d}{2}\frac{p-1}{p+1}\Big]b_0\int_{\R^d} w^{p+1}.\\
\end{split}
\end{equation*}
Since $1-\frac{d}{2}\frac{p-1}{p+1}\neq0$ when $1<p<1+\frac{4}{d}$, one then has $b_0=0$.
Hence, we complete the proof of \eqref{val:b0.0}.

Further, it follows from \eqref{pid:1} that
\begin{equation*}
\sum^d_{i=1}b_i\int_{\R^d}\frac{\partial^2 V_0(x+y_0)}{\partial x_i\partial x_j}w^2=\sum^d_{i=1}b_i\frac{\partial^2 Q(y_0)}{\partial x_i\partial x_j}=0,
\end{equation*}
which implies from the non-degeneracy assumption of $Q(y_0)$ that $b_i=0$ for $i=1,2,...d$.
From \eqref{val:eta0} we thus have
$\bar{\eta}_0\equiv 0$  on $\R^d$.

At last, we claim that $\bar{\eta}_0=0$ cannot occur.
Suppose $\bar{y}_k\in\R^d$ is a maximum point of $\bar{\eta}_k$, and then $|\bar{\eta}_k(\bar{y}_k)|=\|\bar{\eta}_k\|_{L^\infty(\R^d)}=1$.
It thus follows from \eqref{equ:bar.etak} that $\bar{g}_k(\bar{y}_k)+\bar{f}_k(\bar{y}_k)\geq\frac{1}{2}$.
One can  further deduce from \eqref{lim:gk} and \eqref{lim:fk} that $w(\bar{y}_k)\geq C_0>0$, which implies that $y_k$ is bounded uniformly in $k$, due to the fact that $w(x)$ decays exponentially as $|x|\to\infty$.
Therefore, one can  conclude from \eqref{lim:etak} that $\bar{\eta}_0\not\equiv 0$ on $\R^d$, which however contradicts to the fact that $\bar{\eta}_0 \equiv 0$ on $\R^d$. Therefore, the proof of Theorem \ref{Thm:uniqueness} is complete.
\qed

\vskip 0.2truein

\appendix
\section{Appendix}
\subsection{Equivalence between ground states and constraint minimizers}  \label{sec:g.m}

This section is devoted to the proof of Theorem \ref{Thm:g.m} on the equivalence between ground states of equation \eqref{equ:u} and constraint minimizers of problem \eqref{def:e}. At first, we give the  following  lemma. 

\begin{lem}\label{Lem:proper.ea}
Suppose $V(x)$ satisfies \eqref{con:V1}, and $u_\rho(x)\geq0$ is a nonnegative minimizer of $e(\rho)$.
For any $\rho_1,\rho_{2k}\in(0,\infty)$ satisfying $\rho_{2k}\to \rho_1$ as $k\to\infty$, passing to a subsequence if necessary, there exists $ \bar u \in M_{\rho_1} $ such that
\begin{equation}\label{lim:ua1.ua2}
  u_{\rho_{2k}}\to \bar u \quad \text{in}\,\,\, \h(\R^d)\,\,\,\text{as}\,\,\, k\to\infty.
\end{equation}
\end{lem}

\noindent\textbf{Proof.}
For any $\rho_1,\rho_{2k}\in(0,\infty)$, we have
\[ E_{\rho_1}(u_{\rho_1})-E_{\rho_{2k}}(u_{\rho_1})\leq e(\rho_1)-e(\rho_{2k})\leq E_{\rho_1}(u_{\rho_{2k}})-E_{\rho_{2k}}(u_{\rho_{2k}}).\]
One can thus derive that
\begin{equation}\label{ineq:ea1.ea2}
-\frac{\rho_{2k}^{p-1}-\rho_1^{p-1}}{p+1}\int_{\R^d}|u_{\rho_{2k}}|^{p+1}\dx
     \leq  e(\rho_{2k})-e(\rho_1)
     \leq- \frac{\rho_{2k}^{p-1}-\rho_1^{p-1}}{p+1}\int_{\R^d}|u_{\rho_1}|^{p+1}\dx,
\end{equation}
which implies that $e(\rho)$ is a decreasing function of  $\rho\in(0,\infty)$
and $\lim\limits_{k\to\infty}e(\rho_{2k})=e(\rho_1)$,
because  $\|u_{\rho_{2k}}\|_{p+1}^{p+1}$ is bounded uniformly in $k$.

From \eqref{def:E}, one can further deduce that
\begin{equation}\label{lim:Ea1,ua2}
  E_{\rho_1}(u_{\rho_{2k}})=e(\rho_{2k}) + \frac{\rho_{2k}^{p-1}-\rho_1^{p-1}}{p+1} \int_{\R^d} |u_{\rho_{2k}}|^{p+1} dx\to e(\rho_1)\,\,\, \text{as} \,\,\, k\to\infty.
\end{equation}
Let $\{u_{\rho_{2k}}\}$ be a minimizing sequence for $e(\rho_1)$. Employing the  Gagliardo-Nirenberg inequality \eqref{ineq:GNw}, one can derive from \eqref{def:E} that $\{u_{\rho_{2k}}\}$ is bounded uniformly in $\h(\R^d)$ with respect to $k$.
Applying the well-known compact embedding theorem (cf.  \cite[Theorem XIII.67]{RS}), one can deduce that passing to subsequence if necessary, $u_{\rho_{2k}}\to\bar u$ strongly in $L^q(\R^d)$ with $q\in[2,2^*)$ for some $\bar u\in\h$.
This gives the weak lower-semicontinuity  of $E_{\rho_1}(u_{\rho_{2k}})$, and implies from \eqref{lim:Ea1,ua2} that
\[e(\rho_1)=\lim_{k\to\infty}E_{\rho_1}(u_{\rho_{2k}})\geq E_{\rho_1}(\bar u)\geq e(\rho_1), \]
$i.e.$,  (\ref{lim:ua1.ua2}) holds. Hence, the proof of this lemma is completed.
\qed

Next, we giving the following lemma on the differentiability of $e(\rho)$.
\begin{lem}\label{Lem:diff.ea}
Suppose $V(x)$ satisfies \eqref{con:V1}, and let $u_\rho(x)\geq0$ be a nonnegative minimizer of $e(\rho)$.   Then
$e(\rho)$ is differentiable for $a.e.$ $\rho\in(0,\infty)$  and
\begin{equation}\label{diff.ea}
 e'(\rho)=-\frac{p-1}{p+1}\rho^{p-2}\int_{\R^d}|u|^{p+1} \dx,\,\,\  \forall\ u \in M_\rho.
\end{equation}
\end{lem}

Since the proof of this lemma is similar to that of \cite[Lemma 2.2]{GWZZ}, we omit it here.
\qed

Based on the proof of lemma \ref{Lem:diff.ea}, we remark that for any given $\rho\in(0,\infty)$, if $e(\rho)$ admits a unique nonnegative or nonpositive minimizer, then $ e'(\rho)$ exists and satisfies \eqref{diff.ea}.

\vskip 0.2truein

\noindent\textbf{Proof of Theorem  \ref{Thm:g.m}:}
For any $\rho\in(0,\infty)$ and $0\leq u_\rho\in M_\rho$,  $ u_\rho $ satisfies  (\ref{equ:u}) for some Lagrange multiplier $\mu_\rho\in\R$.
It then follows from \eqref{equ:u}, \eqref{def:e} and (\ref{diff.ea}) that, for $a.e.\ \rho\in(0,\infty)$,
\begin{equation}\label{mu.e}
\mu_\rho=2e(\rho)-\frac{p-1}{p+1}\rho^{p-1}\int_{\R^d}|u|^{p+1} \dx=2e(\rho)+\rho e'(\rho),
\end{equation}
which implies that $\mu_\rho$ depends only on the value of $\rho$ and is  independent of the choice of $u_\rho$.
This further indicates that, for $a.e.\ \rho \in (0, \infty)$, all minimizers of $e(\rho)$ satisfy equation (\ref{equ:u}) with the same Lagrange multiplier $\mu_\rho$.

Taking any  $u_g\in G_{\rho,\mu}$ and setting $\tilde{u}_g=\frac{1}{\|u_g\|_2}u_g$,
one then has $F_{\rho,\mu}(\tilde{u}_g) \geq F_{\rho,\mu}(u_g)$.
Since $u_g$ solves \eqref{equ:u}, one can derive from \eqref{def:F} that
\begin{equation*}
F_{\rho,\mu}(u_g)=\Big(\frac{1}{2}-\frac{1}{p+1}\Big)\rho^{p+1}\int_{\R^d}|u_g|^{p+1} \dx
\end{equation*}
and
\begin{equation*}
  F_{\rho,\mu}(\tilde{u}_g)
=\Big(\frac{1}{2\|u_g\|_2^2}-\frac{1}{(p+1)\|u_g\|_2^{p+1}}\Big) \rho^{p+1}\int_{\R^d}|u_g|^{p+1} \dx.
\end{equation*}
Therefore, we have
\begin{equation}\label{ff}
\frac{1}{2\|u_g\|_2^2}-\frac{1}{(p+1)\|u_g\|_2^{p+1}}\geq\frac{1}{2}-\frac{1}{p+1}.
\end{equation}
One can check that \eqref{ff} holds if and only if $\|u_g\|_2=1$, $i.e.$, $F_{\rho,\mu}(\tilde{u}_g)=F_{\rho,\mu}(u_g)$.
This implies that all ground states of equation \eqref{equ:u} share the same $L^2$-norm, $i.e.$,
\[\text{for any $u_g\in G_{\rho,\mu}$, $u_g$ satisfies $\|u_g\|_2^2=1$.}\]

For any $u_g\in G_{\rho,\mu}$ and $u_\rho\in M_\rho$, one has
\begin{equation*}
\text{$E_\rho(u_g)\geq E_\rho(u_\rho)$ and $F_{\rho,\mu}(u_\rho) \geq F_{\rho,\mu}(u_g)$.}
\end{equation*}
Following from \eqref{def:E} and \eqref{def:F}, one has
\begin{equation}\label{FE}
\text{$F_{\rho,\mu}(u)=E_\rho(u)-\frac{1}{2}\mu$,}
\end{equation}
which indicates that $E_{\rho}(u_\rho) \geq E_{\rho}(u_g)$, $i.e.$, $u_g\in M_\rho$.
One can further deduce from \eqref{mu.e} that for $a.e.\,\rho\in(0,\infty)$, there holds that $\mu=\mu_\rho$, which implies $F_{\rho,\mu_\rho}(u_g) \geq F_{\rho,\mu_\rho}(u_\rho)$, $i.e.$, $u_\rho\in G_{\rho,\mu_\rho}$.
Hence we complete the proof of Theorem \ref{Thm:g.m}.
\qed

\vskip 0.2truein

\subsection{Some results on the problem $\tilde{e}_\rho$}\label{sec:tile}
In this section, we focus on studying the following minimization problem
\begin{equation}\label{def:tile}
  \tilde{e}(\rho):=\inf\big\{\tilde{E}_{\rho}(u):u\in H^1(\R^d),\|u\|_{2}=1\big\},
\end{equation}
where $\tilde{E}_{\rho}(u)$ is defined by
\begin{equation}\label{def:tilE}
 \tilde{E}_{\rho}(u):=\frac{1}{2}\int_{\mathbb{R}^d}|\nabla u|^2\dx
  -\frac{\rho^{p-1}}{p+1}\int_{\mathbb{R}^d}|u|^{p+1}\dx, \,\,\,1<p<1+\frac{4}{d}.
\end{equation}
Employing the concentration-compactness principle, one can derive that $\tilde{e}(\rho)$ admits minimizers for any $\rho\in(0,\infty)$, see, e.g., \cite{C,Lions1,Lions2}.
Similar to problem \eqref{def:e}, without loss of generality, we restrict the minimizers of problem \eqref{def:tile} to nonnegative functions.
We then give our results by the following lemma.

\begin{lem}\label{iem:tile}
Suppose $\tilde{u}_\rho$ is a nonnegative minimizer of $\tilde{e}(\rho)$.
Set  $\tilde{\alpha}_\rho:=\Big(\frac{\rho}{\sqrt{a^*}}\Big)^{\frac{2(p-1)}{4-d(p-1)}}$ and
$a^*:=\|w\|_2^2$, where  $w$ is the unique positive solution of \eqref{equ:w}. We then have
\begin{equation}\label{val:tile}
\tilde{e}(\rho)=-\lambda \Big(\frac{\rho}{\sqrt{a^*}}\Big)^\frac{4(p-1)}{4-d(p-1)},
\end{equation}
and, up to translations, $\tilde{u}_\rho$ satisfies
\begin{equation}\label{val:tilu}
\tilde{u}_\rho:=\frac{1}{\sqrt{a^*}}\tilde{\alpha}_\rho^\frac{d}{2}w(\alpha_\rho x),
\end{equation}
where $\lambda$ is defined in \eqref{def:lam}.
\end{lem}
\noindent {\bf Proof.} Suppose $\tilde{u}_\rho$ is a nonnegative minimizer of $\tilde{e}(\rho)$ and $\tilde{u}_1$ is a nonnegative minimizer of $\tilde{e}(1)$.
At first,  we claim that
\begin{equation}\label{eq:erho.e1}
\tilde{e}(\rho)=\rho^\frac{4(p-1)}{4-d(p-1)}\tilde{e}(1)
\,\,\,\text{and}\,\,\,
\tilde{u}_\rho=\alpha_\rho^\frac{d}{2}\tilde{u}_1(\alpha_\rho x),
\end{equation}
 where $\alpha_\rho:=\rho^{\frac{2(p-1)}{4-d(p-1)}}$.
In fact, setting $\tilde{v}_1:=\alpha_\rho^{-\frac{d}{2}}\tilde{u}_\rho(\alpha_\rho^{-1} x)$,
one can check that
\begin{equation*}
  \tilde{e}(\rho)= \tilde{E}_{\rho}(\tilde{u}_\rho)
=\rho^{\frac{4(p-1)}{4-d(p-1)}}\Big[\frac{1}{2}\int_{\mathbb{R}^d}|\nabla \tilde{v}_1|^2\dx
  -\frac{1}{p+1}\int_{\mathbb{R}^d}\tilde{v}_1^{p+1}\dx\Big]
\geq \rho^{\frac{4(p-1)}{4-d(p-1)}}\tilde{e}(1).  \\
\end{equation*}
Similarly, setting $\tilde{v}_\rho:=\alpha_\rho^\frac{d}{2}\tilde{u}_1(\alpha_\rho x)$,  one can check that
\begin{equation}
  \tilde{e}(\rho)\leq\tilde{E}_{\rho}(\tilde{v}_\rho)
=\rho^{\frac{4(p-1)}{4-d(p-1)}}\tilde{e}(1).  \\
\end{equation}
The above two inequalities then give the first equality in \eqref{eq:erho.e1}.
Furthermore, we know that $\tilde{v}_1$ is a minimizer of $\tilde{e}(1)$ and $\tilde{v}_\rho$ is a minimizer of $\tilde{e}(\rho)$, which gives the second equality in \eqref{eq:erho.e1}.

Next, we claim that
\begin{equation}\label{val:tile1}
\tilde{e}(1)=-\lambda(\sqrt{a^*})^{-\frac{4(p-1)}{4-d(p-1)}},
\end{equation}
where $\lambda$ is given by \eqref{def:lam}, and $\tilde{u}_\rho$ (up to translations) satisfies
\begin{equation}\label{val:tilu1}
\tilde{u}_1(x)=(\sqrt{a^*})^{-\frac{4}{4-d(p-1)}}w\big((\sqrt{a^*})^{-\frac{2(p-1)}{4-d(p-1)}}x\big).
\end{equation}
Take a test function $\tilde{v}_\epsilon=\epsilon^\frac{d}{2}\tilde{v}_0(\epsilon x)$, where $0<\tilde{v}_0\in H^1(\R^d)$ satisfies $\|\tilde{v}_0\|_2^2=1$ and $\epsilon>0$ is a positive constant.  One can then verify that
\begin{equation*}
\tilde{e}(1)\leq\tilde{E}_{1}(\tilde{v}_\epsilon)=
\frac{\epsilon^2}{2}\int_{\mathbb{R}^d}|\nabla \tilde{v}_0|^2\dx
  -\frac{\epsilon^{\frac{d}{2}(p-1)}}{p+1}\int_{\mathbb{R}^d}|\tilde{v}_0|^{p+1}\dx<0,\,\,\,\text{when $\epsilon$ is small enough.}
\end{equation*}
Let $\tilde{u}_1$ be a nonnegative minimizer of $\tilde{e}(1)$, and then $\tilde{u}_1$ solves
\begin{equation}\label{equ:tilu1}
  -\Delta \tilde{u}_1=\tilde{\mu}_1\tilde{u}_1+\tilde{u}_1^p,
\end{equation}
where $\tilde{\mu}_1$ is a suitable Lagrange parameter.
It follows from \eqref{def:tilE} and \eqref{equ:tilu1} that
\begin{equation}\label{val:tilmu1.e}
\tilde{\mu}_1=2\tilde{e}(1)-\frac{p-1}{p+1}\int_{\mathbb{R}^d}\tilde{u}_1^{p+1}\dx<0.
\end{equation}
Applying the maximum principle (cf. \cite{GT}) then yields that $\tilde{u}_1>0$,
which implies that, up to translations,
\[\tilde{u}_1=(-\tilde{\mu}_1)^\frac{1}{p-1}w\big((-\tilde{\mu}_1)^{\frac{1}{2}} x\big),\]
due to the fact that $w$ is the unique positive solution of \eqref{equ:w}.
Furthermore, since $\|\tilde{u}_1\|_2^2=1$, one can then deduce that $\mu_1$ satisfies
\begin{equation*}\label{val:tilmu1.w}
\mu_1=-\|w\|_2^{-\frac{4(p-1)}{4-d(p-1)}}=-(\sqrt{a^*})^{-\frac{4(p-1)}{4-d(p-1)}},
\end{equation*}
which implies \eqref{val:tilu1}.
Further,  substituting \eqref{val:tilu1} into \eqref{def:tilE} then yields \eqref{val:tile1}.

Combining the above two claims then yields \eqref{val:tile} and \eqref{val:tilu}, and this completes the proof of Lemma \ref{iem:tile}.
\qed

\vskip 0.3truein

\noindent {\bf Acknowledgements:} The authors are grateful to Yujin Guo for his fruitful discussions on the present paper.

\end{document}